\documentclass[a4paper,11pt]{article}
\usepackage[top=80pt, bottom=90pt, left=70pt, right=70pt]{geometry}

\usepackage{amsmath,amssymb}
\usepackage{mathtools}
\usepackage{amsthm,thmtools}
\usepackage{enumitem}
\usepackage{algorithm}
\usepackage{algpseudocode}
\algrenewcommand\algorithmicrequire{\textbf{Input:}}
\algrenewcommand\algorithmicensure{\textbf{Output:}}
\usepackage[%
  bookmarksnumbered,
  hypertexnames=false,
  pdfdisplaydoctitle,
  pdfusetitle,
  unicode,
  colorlinks=true,
  linkcolor=blue,
  citecolor=blue
]{hyperref}
\usepackage[capitalise,noabbrev]{cleveref}
\usepackage{bookmark}

\usepackage{mleftright}
\mleftright

\declaretheorem[style=plain,numberwithin=section,name=Theorem]{theorem}
\declaretheorem[style=plain,sibling=theorem,name=Lemma]{lemma}
\declaretheorem[style=plain,sibling=theorem,name=Proposition]{proposition}
\declaretheorem[style=plain,sibling=theorem,name=Corollary]{corollary}
\declaretheorem[style=plain,sibling=theorem,name=Claim]{claim}

\declaretheorem[style=definition,sibling=theorem,name=Remark,qed=$\blacksquare$]{remark}

\crefname{theorem}{Theorem}{Theorems}
\crefname{lemma}{Lemma}{Lemmas}
\crefname{proposition}{Proposition}{Propositions}
\crefname{corollary}{Corollary}{Corollaries}
\crefname{claim}{Claim}{Claims}
\crefname{definition}{Definition}{Definitions}
\crefname{example}{Example}{Examples}
\crefname{remark}{Remark}{Remarks}
\crefname{section}{Section}{Sections}
\crefname{algorithm}{Algorithm}{Algorithms}

\numberwithin{equation}{section}

\newcommand{\Z}{\mathbb{Z}}
\newcommand{\Q}{\mathbb{Q}}
\newcommand{\R}{\mathbb{R}}
\newcommand{\F}{\mathbb{F}}
\newcommand{\K}{\mathbb{K}}
\newcommand{\E}{\mathbb{E}}

\DeclareMathOperator{\supp}{supp}

\DeclarePairedDelimiter{\abs}{\lvert}{\rvert}
\DeclarePairedDelimiter{\rbra}{\lparen}{\rparen}
\DeclarePairedDelimiter{\cbra}{\lbrace}{\rbrace}
\DeclarePairedDelimiter{\ceil}{\lceil}{\rceil}
\DeclarePairedDelimiterX{\Set}[2]{\lbrace}{\rbrace}{#1\;\delimsize\vert\;#2}

\DeclareMathOperator{\rank}{rank}
\DeclareMathOperator{\ncrk}{nc\text{-}rank}
\DeclareMathOperator{\im}{im}
\DeclareMathOperator{\spn}{span}
\DeclareMathOperator{\End}{End}
\DeclareMathOperator{\GL}{GL}
\DeclareMathOperator{\Det}{Det}
\newcommand{\D}{\mathcal{D}}
\newcommand{\Dop}{\D^{\mathrm{op}}}
\newcommand{\op}{\mathrm{op}}
\newcommand{\Mat}{M}
\newcommand{\bX}{\mathbf{X}}
\newcommand{\calM}{\mathcal{M}}
\newcommand{\calA}{\mathcal{A}}
\newcommand{\ratfield}[1]{#1(x_0, x_1, \dots, x_m)}
\newcommand{\freefield}[1]{#1(\langle x_0, x_1, \dots, x_m \rangle)}
\newcommand{\freealg}[1]{#1\langle x_0, x_1, \dots, x_m \rangle}
\newcommand{\spl}{\Phi_{\mathrm{sp}}}
\newcommand{\emb}{\iota}
\newcommand{\bT}{\mathbf{T}}
\newcommand{\bitsize}[1]{\langle #1 \rangle}
\newcommand{\calP}{\mathcal{P}}
\newcommand{\calQ}{\mathcal{Q}}
\newcommand{\calO}{\mathcal{O}}
\DeclareMathOperator{\conv}{conv}
\DeclareMathOperator{\proj}{proj}

\allowdisplaybreaks

\begin{document}

\title{Blow-up size for the noncommutative rank\\
of linear symbolic matrices with low-rank coefficients}
\author{Yuni Iwamasa\thanks{Department of Mathematics, Graduate School of Science, Kobe University, Kobe, 657-8501, Japan.
    Email: \texttt{iwamasa@math.kobe-u.ac.jp}}}
\date{\today}
\maketitle

\begin{abstract}
  The noncommutative rank (nc-rank) of a linear symbolic matrix $A = A_0 x_0 + A_1 x_1 + \dots + A_m x_m$
  is the supremum of $1/k$ times the rank of the $k$-th blow-up $\sum_i A_i \otimes X_i$ of $A$, where $X_i$ are $k \times k$ matrices consisting of distinct variables.
  In this paper, we show that the blow-up size needed to attain the nc-rank is bounded by the ranks of the coefficient matrices.
  More precisely, if $A_1, \dots, A_m$ have rank at most $r$, with no restriction on $A_0$,
  then the $k$-th blow-up attains the nc-rank for every $k \geq r$.
  This holds for rectangular matrices over an arbitrary field,
  and unifies and extends the known results that $k = 1$ suffices for rank-one coefficients (Lov\'asz, 1989; Ivanyos, Karpinski, and Saxena, 2010)
  and that $k = 2$ suffices for skew-symmetric coefficients of rank two (Oki and Soma, 2025).
  We present two applications.
  First, if $r$ is a constant, then the nc-rank of an $n \times n$ linear symbolic matrix can be computed by a randomized algorithm
  with $O(m n^2 + n^\omega)$ arithmetic operations over a sufficiently large field, where $\omega$ is the matrix multiplication exponent.
  Second, the hierarchy of polytopes introduced by Hirai (2026), which converges to the nc-independent set polytope, stabilizes at the $r$-th level;
  this is based on a counterpart of our main result for the degree of the Dieudonn\'e determinant.
\end{abstract}
\begin{quote}
  {\bfseries Keywords:}
  noncommutative rank, linear symbolic matrix, blow-up, randomized algorithm, Dieudonn\'e determinant.
\end{quote}

\section{Introduction}\label{sec:intro}
Let $\F$ be a field and let
\begin{align}\label{eq:linear-matrix}
  A \coloneqq A_0 x_0 + A_1 x_1 + \dots + A_m x_m
\end{align}
be a \emph{linear symbolic matrix} over $\F$,
where $A_0, A_1, \dots, A_m \in \F^{a \times b}$ and $x_0, x_1, \dots, x_m$ are variables.
If the variables commute, then the rank of $A$ is considered over the rational function field $\ratfield{\F}$,
and computing it is known as \emph{Edmonds' problem}~\cite{Edmonds1967}.
Edmonds' problem includes fundamental problems in combinatorial optimization,
such as the bipartite matching, linear matroid intersection, and linear matroid parity problems, as special cases~\cite{Edmonds1967,Lovasz1989}.
It is solved in randomized polynomial time by substituting random values for the variables if $\F$ is sufficiently large~\cite{Lovasz1989},
whereas the existence of a deterministic polynomial-time algorithm is a major open problem in theoretical computer science, which is related to circuit lower bounds~\cite[Theorem~4.12]{Kabanets2004};
see also~\cite{Garg2020,Ivanyos2017}.
If the variables are pairwise noncommutative, then the rank of $A$ is considered over the free skew field $\freefield{\F}$~\cite{Amitsur1966,Cohn1995},
and is called the \emph{noncommutative rank} (\emph{nc-rank}) of $A$, denoted by $\ncrk A$;
computing it is known as the \emph{noncommutative Edmonds' problem}~\cite{Ivanyos2017}.
The nc-rank is an upper bound on the commutative rank,
and, in contrast to the commutative rank, it is known to be computable in deterministic polynomial time~\cite{Garg2020,Hamada2021,Ivanyos2018}.
The noncommutative Edmonds' problem arises in various areas such as invariant theory, quantum information theory, and optimization~\cite{Derksen2017,Garg2020,Gurvits2004,Hamada2021,Ivanyos2017},
and has led to the framework of \emph{algebraic combinatorial optimization}~\cite{Hirai2025}.

The nc-rank admits a characterization via \emph{blow-ups}.
For a positive integer $k$, the \emph{$k$-th blow-up} of $A$ is the matrix
\begin{align*}
  A^{(k)} \coloneqq \sum_{i=0}^m A_i \otimes X_i^{(k)},
\end{align*}
where $X_0^{(k)}, X_1^{(k)}, \dots, X_m^{(k)}$ are $k \times k$ matrices whose entries are distinct commuting variables.
Then $\ncrk A$ is equal to the supremum of $(\rank A^{(k)})/k$ over all $k \geq 1$~\cite{Ivanyos2017,Hirai2025},
where the rank of $A^{(k)}$ is taken over the rational function field, denoted by $\F(\bX^{(k)})$, in the entries of $X_0^{(k)}, \dots, X_m^{(k)}$.

A natural question is how large $k$ must be so that the $k$-th blow-up \emph{attains} the nc-rank, that is, $\ncrk A = (\rank A^{(k)})/k$.
Derksen and Makam~\cite[Theorem~1.8]{Derksen2017} showed that, for $n \times n$ linear symbolic matrices, every $k \geq n - 1$ suffices (see also~\cite[Corollary~3.2]{Derksen2018}).
For special classes of coefficient matrices, much smaller blow-ups are sufficient.
If all of $A_1, \dots, A_m$ have rank one (and $A_0$ is arbitrary), then $k = 1$ suffices,
that is, the commutative rank and the nc-rank coincide~\cite{Lovasz1989,Ivanyos2010}.
If all coefficients are skew-symmetric matrices of rank two, then $k = 2$ suffices~\cite[Theorem~4.8]{Oki2025}.

In these two special cases, the blow-up size attaining the nc-rank is bounded by the rank of the coefficient matrices,
independently of the size of the linear symbolic matrix.
In this paper, we show that such a bound on the blow-up size in terms of the ranks of the coefficient matrices holds in general.

\begin{theorem}\label{thm:main}
  Let $\F$ be a field, and $r$ a positive integer.
  Let $A = A_0 x_0 + \sum_{i=1}^m A_i x_i$ be a linear symbolic matrix with $A_0, A_1, \dots, A_m \in \F^{a \times b}$ such that $\rank A_i \leq r$ for all $i = 1, \dots, m$;
  the matrix $A_0$ is arbitrary.
  Then, for every integer $k \geq r$, we have
  \begin{align}\label{eq:main}
    \ncrk A = \frac{1}{k} \rank A^{(k)}.
  \end{align}
\end{theorem}

In particular, the $r$-th blow-up attains the nc-rank.
Setting $A_0 = 0$ yields the case where all coefficients have rank at most $r$,
and setting $r = k = 1$ recovers the coincidence of the commutative rank and the nc-rank for rank-one coefficients with an arbitrary constant term~\cite{Ivanyos2010}.
Note that the rank on the right-hand side of~\eqref{eq:main} is the generic rank over $\F(\bX^{(k)})$;
if $\F$ is a small finite field, it is not necessarily attained by substituting matrices over $\F$ for $X_0^{(k)}, X_1^{(k)}, \dots, X_m^{(k)}$.

The \emph{matrix space} spanned by the coefficients of $A$ is the $\F$-subspace of $\F^{a \times b}$ spanned by $A_0, A_1, \dots, A_m$.
It is known~\cite[Corollary~1 and Lemma~1]{Fortin2004} that the rank of each blow-up, and hence the nc-rank, depends only on this matrix space; see also \Cref{lem:matrix-space}.
Therefore, it suffices that the rank condition in \Cref{thm:main} is satisfied by some matrices spanning the same matrix space.

\begin{corollary}\label{cor:matrix-space}
  Let $\F$ be a field, and $r$ a positive integer.
  Let $A$ be a linear symbolic matrix over $\F$ of the form~\eqref{eq:linear-matrix}.
  Suppose that the matrix space spanned by $A_0, A_1, \dots, A_m$ is also spanned by some matrices $B_0, B_1, \dots, B_{m'} \in \F^{a \times b}$
  such that $\rank B_j \leq r$ for all $j = 1, \dots, m'$;
  the matrix $B_0$ is arbitrary.
  Then~\eqref{eq:main} holds for every integer $k \geq r$.
\end{corollary}

\Cref{thm:main} also yields an alternative proof of the above-mentioned bound of Derksen and Makam, in the following form (see also~\cite[Section~2.1]{Huszar2021}):
For every linear symbolic matrix $A$, with no assumption on the ranks of the coefficients,
the $k$-th blow-up attains the nc-rank for every positive integer $k \geq \ncrk A - 1$,
and in particular for every $k \geq n - 1$ if $A$ is of size $n \times n$ (\Cref{cor:general-bound}).

We present two applications of \Cref{thm:main}.
First, it yields a simple randomized algorithm for computing the nc-rank under the assumption of \Cref{cor:matrix-space},
in which the matrices $B_0, B_1, \dots, B_{m'}$ need not be given.
For an $n \times n$ linear symbolic matrix and a constant $r$,
this algorithm outputs the nc-rank with probability at least $1/2$ in $O(m n^2 + n^\omega)$ arithmetic operations, provided that $\F$ is sufficiently large,
where $2 \leq \omega \leq 3$ is the matrix multiplication exponent (\Cref{thm:randomized});
a randomized algorithm with the same running time was known for the case where the matrix space has a basis consisting of skew-symmetric matrices of rank two~\cite{Oki2025}.
Second, we show in \Cref{thm:newton} that, under the assumption of \Cref{thm:main},
the hierarchy of polytopes introduced by Hirai~\cite{Hirai2026},
in which the $k$-th polytope is defined from the exponent vectors of the monomials appearing in the minors of the $k$-th blow-up,
stabilizes at the $r$-th level;
this extends the result of Hirai~\cite[Theorem~3.13]{Hirai2026} for the case where the coefficients are skew-symmetric matrices of rank two.
The proof relies on \Cref{prop:deg-det-attainment},
which is a counterpart of \Cref{thm:main} for the degree of the \emph{Dieudonn\'e determinant}~\cite{Dieudonne1943,Hirai2019} and may be of independent interest.

\paragraph{Organization.}
The rest of this paper is organized as follows.
In \Cref{sec:preliminaries}, we introduce notation, known characterizations of the nc-rank, and the tools on division algebras that we use.
\Cref{sec:proof} is devoted to the proofs of \Cref{thm:main}, \Cref{cor:matrix-space}, and a consequence for general linear symbolic matrices.
In \Cref{sec:application}, we present two applications of \Cref{thm:main}:
fast randomized computation of the nc-rank and the stabilization of a hierarchy of polytopes converging to the nc-independent set polytope.

\section{Preliminaries}\label{sec:preliminaries}
For a positive integer $m$, let $[m] \coloneqq \cbra{1, 2, \dots, m}$ and $[0, m] \coloneqq \cbra{0, 1, \dots, m}$.
All fields are commutative, and a \emph{skew field} is a (not necessarily commutative) division ring.
For a vector space $W'$ over a field (resp., a right vector space $W'$ over a skew field),
we write $W \leq W'$ to mean that $W$ is a subspace of $W'$ (resp., a right subspace).
We use the conventions that a zero-dimensional space has the empty basis,
and that a square matrix of size $0$ is nonsingular with determinant $1$.

A matrix over a field $\F$ is often regarded also as a matrix over an extension field of $\F$ or over a skew field containing $\F$,
and a vector space over an extension field of $\F$ is also a vector space over $\F$.
We use subscripts to specify over which field or skew field an object is considered.
For a field $\F$, we denote by $\rank_\F B$, $\ker_\F B$, and $\im_\F B$ the rank, the kernel, and the image of a matrix $B$ regarded as a matrix over $\F$, respectively,
by $\dim_\F W$ the dimension of $W$ as an $\F$-vector space,
and by $\spn_\F S$ the $\F$-linear span of a set $S$ of vectors.
We use the same notation for a skew field $\D$, for which the precise definitions are given in \Cref{subsec:skew-fields}.
We omit the subscript if it is clear from the context;
for instance, $\rank B$ for a matrix $B$ over $\F$ means $\rank_\F B$ unless otherwise stated.

\subsection{Noncommutative rank and blow-ups}\label{subsec:nc-rank}
Let $\F$ be a field, and $A$ a linear symbolic matrix over $\F$ of the form~\eqref{eq:linear-matrix},
where $A_0, A_1, \dots, A_m \in \F^{a \times b}$.
We denote by $\ratfield{\F}$ the rational function field over $\F$ in commuting variables $x_0, x_1, \dots, x_m$;
the rank of $A$ over $\ratfield{\F}$ is referred to as the \emph{commutative rank} of $A$.
We denote by $\freealg{\F}$ the free algebra generated by pairwise noncommutative variables $x_0, x_1, \dots, x_m$,
and by $\freefield{\F}$ its free skew field~\cite{Amitsur1966,Cohn1995}.
The \emph{nc-rank} (or \emph{noncommutative rank}) of $A$ is defined as the rank of $A$ over $\freefield{\F}$, i.e.,
\begin{align*}
  \ncrk_\F A \coloneqq \rank_{\freefield{\F}} A.
\end{align*}
We omit the subscript $\F$ if it is clear from the context.
If $a = 0$ or $b = 0$, then all the ranks considered in this paper are $0$ and \Cref{thm:main} trivially holds;
hence we assume $a, b \geq 1$ in the following.

For a positive integer $k$ and $i \in [0, m]$,
let $X_i^{(k)}$ be a $k \times k$ matrix whose entries are variables,
where all the $(m+1)k^2$ variables are distinct and commute with each other.
We denote by $\bX^{(k)} \coloneqq (X_0^{(k)}, X_1^{(k)}, \dots, X_m^{(k)})$ the tuple of these matrices,
by $\F[\bX^{(k)}]$ the polynomial ring over $\F$ in the entries of the matrices in $\bX^{(k)}$,
and by $\F(\bX^{(k)})$ its field of fractions.
The \emph{$k$-th blow-up} of $A$ is the matrix
\begin{align*}
  A^{(k)} \coloneqq \sum_{i=0}^m A_i \otimes X_i^{(k)} \in \F[\bX^{(k)}]^{ak \times bk},
\end{align*}
that is, the $a \times b$ block matrix whose $(j, \ell)$-th block is the $k \times k$ matrix $\sum_{i=0}^m (A_i)_{j\ell} X_i^{(k)}$.
For $k = 1$, by identifying the $1 \times 1$ matrix $X_i^{(1)}$ with the variable $x_i$, we have $A^{(1)} = A$ and $\F(\bX^{(1)}) = \ratfield{\F}$;
hence the rank of the first blow-up is just the commutative rank of $A$.
When we consider an extension field $\F'$ of $\F$,
we always assume that the entries of the matrices in $\bX^{(k)}$ are algebraically independent over $\F'$ as well,
and use the same notation $\F'[\bX^{(k)}]$ and $\F'(\bX^{(k)})$.

We use the following two known characterizations of the nc-rank.

\begin{theorem}[{\cite{Ivanyos2017}; see also \cite[Proposition~7.3]{Hrubes2015} and \cite[Theorem~2.1]{Hirai2025}}]\label{thm:blow-up}
  For a linear symbolic matrix $A$ over a field $\F$, we have
  \begin{align*}
    \ncrk_\F A = \sup_{k \geq 1} \frac{1}{k} \rank_{\F(\bX^{(k)})} A^{(k)}.
  \end{align*}
\end{theorem}

\begin{theorem}[{\cite[Theorem~1]{Fortin2004}}]\label{thm:fortin-reutenauer}
  For a linear symbolic matrix $A$ over a field $\F$ of the form~\eqref{eq:linear-matrix}, we have
  \begin{align*}
    \ncrk_\F A = \min \cbra{b - \dim_\F U + \dim_\F V},
  \end{align*}
  where the minimum is taken over all pairs $(U, V)$ of subspaces $U \leq \F^b$ and $V \leq \F^a$ such that $A_i U \leq V$ for all $i \in [0, m]$.
\end{theorem}

\subsection{Linear algebra over skew fields}\label{subsec:skew-fields}
Let $\D$ be a skew field.
The \emph{center} of $\D$ is $Z(\D) \coloneqq \Set{z \in \D}{pz = zp \text{ for all } p \in \D}$,
which is a subfield of $\D$.
For a nonnegative integer $h$, we regard $\D^h$ as the right $\D$-vector space of column vectors.
For $w \in \D^h$ and $p \in \D$, we denote by $wp$ and $pw$ the vectors obtained from $w$ by multiplying each entry by $p$ from the right and from the left, respectively.
A matrix $B \in \D^{h_1 \times h_2}$ is regarded as the right $\D$-linear map $\D^{h_2} \to \D^{h_1}$ defined by $w \mapsto Bw$.
We denote by $\ker_\D B$ and $\im_\D B$ its kernel and image, respectively,
and define $\rank_\D B \coloneqq \dim_\D \im_\D B$,
where $\dim_\D$ denotes the dimension as a right $\D$-vector space.
Let $\GL_h(\D)$ be the group of nonsingular $h \times h$ matrices over $\D$,
where a square matrix is said to be \emph{nonsingular} if it has an inverse.
For a subset of $\D^h$, let $\spn_\D$ denote its right $\D$-linear span.

A right $\D$-linear map does not preserve the left multiplication in general.
However, if all entries of $C \in \D^{h_1 \times h_2}$ belong to $Z(\D)$, then $C(pw) = p(Cw)$ for all $p \in \D$ and $w \in \D^{h_2}$;
we use this fact repeatedly.

Let $\K$ be a subfield of $Z(\D)$.
Since every element of $\K$ commutes with all elements of $\D$, the left and right multiplications by an element of $\K$ coincide.
Hence $\D$ is a vector space over $\K$, and its multiplication is $\K$-bilinear;
that is, $\D$ is a $\K$-algebra.
More generally, $\D^h$ is also regarded as a vector space over $\K$;
we denote by $\dim_\K$ the dimension as a $\K$-vector space,
and say that a map is $\K$-linear if it is linear with respect to this structure.

The following lemma compares linear algebra over $\K$ and over $\D$.
It consists of standard facts;
we include a proof to make the paper self-contained.

\begin{lemma}\label{lem:scalar-extension}
  Let $\D$ be a skew field, $\K \subseteq Z(\D)$ a subfield, and $h$ a nonnegative integer.
  \begin{enumerate}[label=\textup{(\arabic*)}]
    \item $\xi_1, \dots, \xi_s \in \K^h$ are linearly independent over $\K$ if and only if they are right linearly independent over $\D$ in $\D^h$.
    \item For $W_0 \leq \K^h$, the map $W_0 \otimes_\K \D \to \spn_\D W_0$ defined by $\xi \otimes p \mapsto \xi p$ is an isomorphism of right $\D$-vector spaces.
          In particular, $\dim_\D \spn_\D W_0 = \dim_\K W_0$.
    \item For a matrix $C \in \K^{h \times s}$, we have $\rank_\D C = \rank_\K C$ and $\ker_\D C = \spn_\D \ker_\K C$.
  \end{enumerate}
\end{lemma}
\begin{proof}
  (1).
  The if part is clear.
  For the only-if part, suppose that $\xi_1, \dots, \xi_s$ are linearly independent over $\K$ and that $\sum_{\ell=1}^s \xi_\ell p_\ell = 0$ with $p_\ell \in \D$.
  Take a $\K$-basis $e_1, \dots, e_{s'}$ of the $\K$-linear span of $p_1, \dots, p_s$ in $\D$,
  and write $p_\ell = \sum_{\alpha=1}^{s'} \kappa_{\ell\alpha} e_\alpha$ with $\kappa_{\ell\alpha} \in \K$.
  Then, we have
  \begin{align*}
    0 = \sum_{\ell=1}^{s} \xi_\ell \rbra*{\sum_{\alpha=1}^{s'}\kappa_{\ell\alpha} e_\alpha} = \sum_{\alpha=1}^{s'} \rbra*{\sum_{\ell=1}^s \xi_\ell \kappa_{\ell\alpha}} e_\alpha.
  \end{align*}
  Each vector in the parentheses of the rightmost expression belongs to $\K^h$.
  Then, since $e_1, \dots, e_{s'}$ are linearly independent over $\K$, we have $\sum_{\ell=1}^s \xi_\ell \kappa_{\ell\alpha} = 0$ for each $\alpha$.
  Moreover,
  by the linear independence of $\xi_1, \dots, \xi_s$ over $\K$, we obtain $\kappa_{\ell\alpha} = 0$ for all $\ell$ and $\alpha$.
  Hence $p_\ell = 0$ for all $\ell$.

  (2).
  Let $\xi_1, \dots, \xi_s$ be a $\K$-basis of $W_0$.
  By the assertion~(1), it is a right $\D$-basis of $\spn_\D W_0$.
  The map is clearly surjective.
  Every element of $W_0 \otimes_\K \D$ is written as $\sum_{\ell} \xi_\ell \otimes p_\ell$, and if its image $\sum_{\ell} \xi_\ell p_\ell$ is zero, then $p_\ell = 0$ for all $\ell$ by the right linear independence of $\xi_1, \dots, \xi_s$ over $\D$.
  Hence the map is injective.

  (3).
  The rank equality follows from the assertion~(1), since the rank is the maximum number of (right) linearly independent columns.
  Clearly $\spn_\D \ker_\K C \leq \ker_\D C$,
  and both sides have right dimension $s - \rank_\K C$ by the assertion~(2) and the rank equality.
  Thus, $\ker_\D C = \spn_\D \ker_\K C$ holds.
\end{proof}

\subsection{Two-sided multiplications and central descent}\label{subsec:center}
In this subsection, we assume that $\D$ is a skew field such that $Z(\D) = \K$ and $N \coloneqq \dim_\K \D < \infty$.
Let $\Dop$ denote the opposite algebra of $\D$;
the element of $\Dop$ corresponding to $z \in \D$ is denoted by $z^{\op}$, so that $y^{\op} z^{\op} = (zy)^{\op}$.
Let $\End_\K(\D)$ denote the $\K$-algebra of $\K$-linear maps from $\D$ to itself.
The following is the special case for skew fields of the standard isomorphism for central simple algebras;
we include a proof to make the paper self-contained.

\begin{lemma}[{see e.g.,~\cite[Proof of Proposition~2.4.8]{Gille2006} and~\cite[Tag 074I]{Stacks}}]\label{lem:two-sided-multiplication}
  The $\K$-linear map $\Phi \colon \D \otimes_\K \Dop \to \End_\K(\D)$ defined by $\Phi(y \otimes z^{\op})(p) \coloneqq ypz$ for $p \in \D$ is an isomorphism.
\end{lemma}
\begin{proof}
  We first show the injectivity.
  Suppose, to the contrary, that $\Phi$ is not injective,
  and take a nonzero element $\sum_{\ell=1}^s y_\ell \otimes z_\ell^{\op}$ of the kernel of $\Phi$ such that the number $s$ of terms is minimum.
  Thus $\sum_{\ell=1}^s y_\ell p z_\ell = 0$ for all $p \in \D$.
  By the minimality of $s$, each of the families $y_1, \dots, y_s$ and $z_1, \dots, z_s$ is linearly independent over $\K$;
  otherwise, we can reduce the number of terms by expressing one element as a $\K$-linear combination of the others.
  In particular, $y_1 \neq 0$.
  By multiplying all $y_\ell$ by $y_1^{-1}$ from the left, we may assume $y_1 = 1$;
  this operation preserves the membership in the kernel, the nonzeroness of the element, and the minimality.
  For any $t, p \in \D$, we have $\sum_{\ell=1}^s y_\ell t p z_\ell = 0$ and $\sum_{\ell=1}^s t y_\ell p z_\ell = 0$.
  Since $y_1 = 1$, subtracting the latter from the former yields $\sum_{\ell=2}^s (y_\ell t - t y_\ell) p z_\ell = 0$.
  Hence, for each $t \in \D$, the kernel contains $\sum_{\ell=2}^s (y_\ell t - t y_\ell) \otimes z_\ell^{\op}$, which has at most $s - 1$ terms and hence must be zero by the minimality of $s$.
  Since $z_2, \dots, z_s$ are linearly independent over $\K$, we have $y_\ell t = t y_\ell$ for all $\ell$ and all $t \in \D$,
  that is, $y_\ell \in Z(\D) = \K$ for all $\ell$.
  Since $\sum_{\ell=1}^s y_\ell p z_\ell = 0$ for $p = 1$, we obtain $\sum_{\ell=1}^s y_\ell z_\ell = 0$ with $y_\ell \in \K$ and $y_1 = 1$,
  which contradicts the linear independence of $z_1, \dots, z_s$ over $\K$.
  Therefore $\Phi$ is injective.

  Since both the domain and the codomain of $\Phi$ have $\K$-dimension $N^2$, the map $\Phi$ is also surjective.
\end{proof}

By \Cref{lem:two-sided-multiplication}, every $\K$-linear map on $\D$ is a finite sum of two-sided multiplications $p \mapsto ypz$.
For a $\K$-linear map $G \colon \D \to \D$, let $\rank_\K G$ denote the $\K$-dimension of its image.
The next lemma relates the number of terms and the dimension of the image.

\begin{lemma}\label{lem:operator-bound}
  Let $G \in \End_\K(\D)$ be a map of the form
  \begin{align*}
    G(p) = \sum_{\ell=1}^s y_\ell p z_\ell \qquad (p \in \D)
  \end{align*}
  for some nonnegative integer $s$ and $y_\ell, z_\ell \in \D$.
  If $N > s \rank_\K G$, then $G = 0$.
\end{lemma}
\begin{proof}
  We show the contrapositive:
  If $G \neq 0$, then $N \leq s \rank_\K G$.
  Suppose that $G \neq 0$.
  It suffices to consider the case where $s$ is the minimum number of terms in such an expression of $G$.
  Then $s \geq 1$, the elements $z_1, \dots, z_s$ are linearly independent over $\K$, and $y_\ell \neq 0$ for all $\ell$,
  since otherwise we can reduce the number of terms.
  Let $z \coloneqq (z_1, z_2, \dots, z_s)^\top \in \D^s$ be the column vector over $\D$.

  We first claim that $\spn_\D \Set{pz}{p \in \D}$ is the entire space $\D^s$.
  Suppose not.
  Then, there is a nonzero row vector $(y_1', y_2', \dots, y_s') \in \D^{1 \times s}$ over $\D$ such that
  $\sum_{\ell=1}^s y_\ell' p z_\ell = 0$ for all $p \in \D$.
  By the injectivity in \Cref{lem:two-sided-multiplication}, we have $\sum_{\ell} y_\ell' \otimes z_\ell^{\op} = 0$.
  Since $z_1, z_2, \dots, z_s$ are linearly independent over $\K$, we obtain $y_\ell' = 0$ for all $\ell$, a contradiction.

  By the claim, there are $p_1, \dots, p_s \in \D$ such that $p_1z, p_2z, \dots, p_sz$ form a right $\D$-basis of $\D^s$,
  that is, the matrix $B \in \D^{s \times s}$ whose $(\ell, j)$-th entry is $p_j z_\ell$ belongs to $\GL_s(\D)$.
  Let $t_1, t_2, \dots, t_s \in \D$ be the entries of the first column of its inverse $B^{-1}$.
  Since $B B^{-1} = I_s$, we have $\sum_{j=1}^s p_j z_\ell t_j = 1$ for $\ell = 1$, and $0$ for $\ell \neq 1$.
  Hence, for all $p \in \D$, we have
  \begin{align*}
    \sum_{j=1}^s G(p p_j) t_j = \sum_{j=1}^s \rbra*{\sum_{\ell = 1}^s y_\ell pp_j z_\ell} t_j
    = \sum_{\ell=1}^s y_\ell p \rbra*{\sum_{j=1}^s p_j z_\ell t_j}
    = y_1 p.
  \end{align*}
  This implies that the $\K$-linear maps $p \mapsto y_1 p$ and $p \mapsto \sum_{j=1}^s G(p p_j) t_j$
  coincide, and hence have the same rank over $\K$.
  Since $y_1 \neq 0$,
  the map $p \mapsto y_1 p$ is bijective,
  which implies that its rank over $\K$ is $N$.
  On the other hand, for each $j$, the map $p \mapsto G(p p_j) t_j$ is the composition of $G$ with the right multiplications by $p_j$ and by $t_j$.
  Hence its rank over $\K$ is at most $\rank_\K G$,
  which implies that the rank of the sum $\sum_{j=1}^s G(p p_j) t_j$ over $\K$ is at most $s \rank_\K G$.
  Therefore, we obtain $N \leq s \rank_\K G$.
\end{proof}

A subspace $W \leq \D^h$ is said to be \emph{two-sided} if $pW \leq W$ for all $p \in \D$.
The next tool states that a two-sided subspace is defined over the center;
this is a consequence of~\cite[Tag 074B]{Stacks}.

\begin{lemma}\label{lem:central-descent}
  Let $W \leq \D^h$ be two-sided and $W_0 \coloneqq W \cap \K^h$.
  Then the map $W_0 \otimes_\K \D \to W$ defined by $\xi \otimes p \mapsto \xi p$ is an isomorphism of right $\D$-vector spaces.
  In particular, $W = \spn_\D W_0$ and $\dim_\K W_0 = \dim_\D W$.
\end{lemma}
\begin{proof}
  We first show that $W = \spn_\D W_0$.
  Since $\spn_\D W_0 \leq W$ is clear, we show the converse inclusion.
  Take a $\K$-basis $e_1, \dots, e_N$ of $\D$,
  and expand $w \in W$ as $w = \sum_{\alpha=1}^N \xi_\alpha e_\alpha$ with $\xi_\alpha \in \K^h$.
  For each $\alpha$, let $e_\alpha^* \colon \D \to \D$ be the $\K$-linear map defined by $e_\alpha^*\rbra*{\sum_{\alpha'} \kappa_{\alpha'} e_{\alpha'}} \coloneqq \kappa_\alpha$ for $\kappa_{\alpha'} \in \K$.
  By the surjectivity in \Cref{lem:two-sided-multiplication},
  there are finitely many elements $y_{\alpha\ell}, z_{\alpha\ell} \in \D$ such that $e_\alpha^*(p) = \sum_{\ell} y_{\alpha\ell} p z_{\alpha\ell}$ for all $p \in \D$.
  By applying $e_\alpha^*$ to each entry of $w$, we obtain
  $\xi_\alpha = \sum_{\ell} y_{\alpha\ell} w z_{\alpha\ell}$.
  Since $W$ is two-sided, each term in the right-hand side belongs to $W$.
  Hence, $\xi_\alpha \in W \cap \K^h = W_0$ holds for all $\alpha$,
  and thus $w \in \spn_\D W_0$.

  The assertion now follows from \Cref{lem:scalar-extension}~(2).
\end{proof}

In the following, for a two-sided subspace $W \leq \D^h$, we write $W_0 \coloneqq W \cap \K^h$ and $W = W_0 \otimes_\K \D$,
where the latter is by the identification in \Cref{lem:central-descent}.

A map on a two-sided subspace that is compatible with the multiplications from both sides is also defined over the center.

\begin{lemma}\label{lem:central-map}
  Let $W \leq \D^h$ be two-sided,
  and let $g \colon W \to \D^{h'}$ be a right $\D$-linear map that preserves the left multiplication, that is,
  $g(t w) = t g(w)$ for all $t \in \D$ and $w \in W$.
  Then there is a matrix $C \in \K^{h' \times h}$ such that $g(w) = C w$ for all $w \in W$.
\end{lemma}
\begin{proof}
  Let $\xi \in W_0$.
  Then, for all $t \in \D$, we have $t g(\xi) = g(t \xi) = g(\xi t) = g(\xi) t$,
  where the first equality follows from the fact that $g$ preserves the left multiplication,
  the second from $t \xi = \xi t$, which holds since $\xi \in W_0 \leq \K^h$,
  and the third from the right $\D$-linearity of $g$.
  Thus every entry of $g(\xi)$ belongs to $Z(\D) = \K$.
  Therefore $g$ maps $W_0$ into $\K^{h'}$,
  and its restriction $g|_{W_0} \colon W_0 \to \K^{h'}$ is $\K$-linear, since $g$ is right $\D$-linear.
  We extend it to a $\K$-linear map $\K^h \to \K^{h'}$ by extending a basis of $W_0$ to that of $\K^h$ and sending the additional basis vectors to zero,
  and let $C \in \K^{h' \times h}$ be its matrix representation.
  Then $C$ also defines a right $\D$-linear map $\D^h \ni w \mapsto C w \in \D^{h'}$,
  which coincides with $g$ on $W_0$ by the construction.
  Since both maps are right $\D$-linear and $W = \spn_\D W_0$ by \Cref{lem:central-descent}, they coincide on $W$.
\end{proof}

\subsection{Division algebras of prescribed degree}\label{subsec:division-algebra}
Our proof uses a skew field that is finite-dimensional over its center, with prescribed dimension.

\begin{proposition}[{\cite[Lemma~3.3 and Proposition~4.2]{Ivanyos2018}}]\label{prop:division-algebra}
  Let $\F$ be a field and $d$ a positive integer.
  Then there exist field extensions $\F \subseteq \K \subseteq \E$ and a skew field $\D$ such that
  $\abs{\K} = \infty$, $Z(\D) = \K$, $\dim_\K \D = d^2$,
  and there is an isomorphism of $\E$-algebras $\spl \colon \D \otimes_\K \E \to \Mat_d(\E)$,
  where $\Mat_d(\E)$ denotes the algebra of $d \times d$ matrices over $\E$.
\end{proposition}

The cited results~\cite[Lemma~3.3 and Proposition~4.2]{Ivanyos2018} provide a constructive proof of \Cref{prop:division-algebra};
since \cite[Lemma~3.3]{Ivanyos2018} assumes that the base field contains a certain root of unity,
we apply it after replacing $\F$ with its extension field obtained by adjoining this root of unity.
For existence alone, one can take $\D$ to be the universal division algebra of degree $d$~\cite[Theorem~2.5]{Derksen2018}
over $\F$ (over a rational function field $\F(t)$ if $\F$ is finite),
which is a division algebra of dimension $d^2$ over its center, and whose center is an infinite field containing $\F$;
since $\D$ is then a central simple algebra over its center,
a field $\E$ with the required isomorphism $\spl$ exists by the splitting-field theorem~\cite[Theorem~2.2.1]{Gille2006}.

Let $\F \subseteq \K \subseteq \E$, $\D$, and $\spl$ be as in \Cref{prop:division-algebra}.
We define the embedding $\emb \colon \D \to \Mat_d(\E)$ by $\emb(p) \coloneqq \spl(p \otimes 1)$ for $p \in \D$.
Since $p \mapsto p \otimes 1$ is an injective homomorphism of $\K$-algebras, so is $\emb$.
Moreover, since $\kappa \otimes 1 = 1 \otimes \kappa$ for $\kappa \in \K$ and $\spl$ is $\E$-linear with $\spl(1 \otimes 1) = I_d$,
we have $\emb(\kappa) = \kappa I_d$ for all $\kappa \in \K$.
For a matrix $B = (B_{j\ell}) \in \D^{h_1 \times h_2}$, we denote by $\emb(B) \in \E^{dh_1 \times dh_2}$ the $h_1 \times h_2$ block matrix whose $(j, \ell)$-th block is the $d \times d$ matrix $\emb(B_{j\ell})$.
Then $\emb(BB') = \emb(B)\emb(B')$ holds whenever the product $BB'$ is defined.
The following rank comparison is described in~\cite[Section~2.1.1]{Hirai2025} for a general embedding of a skew field into a matrix algebra;
we include a proof for our setting.

\begin{lemma}\label{lem:splitting-rank}
  For every $B \in \D^{h_1 \times h_2}$, we have $\rank_\E \emb(B) = d \rank_\D B$.
\end{lemma}
\begin{proof}
  Let $\ell \coloneqq \rank_\D B$.
  There are nonsingular matrices $Y_1 \in \GL_{h_1}(\D)$ and $Y_2 \in \GL_{h_2}(\D)$ such that
  \begin{align*}
    Y_1 B Y_2 =
    \begin{pmatrix}
      I_\ell & 0 \\
      0      & 0
    \end{pmatrix}.
  \end{align*}
  Since $\emb$ preserves products and the identity, $\emb(Y_1)$ and $\emb(Y_2)$ are nonsingular over $\E$, and
  \begin{align*}
    \emb(Y_1) \emb(B) \emb(Y_2) =
    \emb\rbra*{Y_1 B Y_2}
    =
    \begin{pmatrix}
      I_{d\ell} & 0 \\
      0         & 0
    \end{pmatrix}.
  \end{align*}
  Hence, we obtain $\rank_\E \emb(B) = d\ell$.
\end{proof}

We also use the following principle:
If perturbing a matrix in a direction never increases the rank, then the direction maps the kernel into the image.
The corresponding statement over a field is~\cite[Lemma~2.2]{Ivanyos2010}.
We prove it for the skew field $\D$ by passing to the splitting field $\E$, while the perturbation parameter is taken from the center $\K$;
this is where the infiniteness of $\K$ is used.

\begin{lemma}\label{lem:max-rank}
  Let $B, B' \in \D^{h_1 \times h_2}$.
  If $\rank_\D (B + \lambda B') \leq \rank_\D B$ for all $\lambda \in \K$, then $B' \ker_\D B \leq \im_\D B$.
\end{lemma}
\begin{proof}
  Let $\ell \coloneqq \rank_\D B$.
  Take $Y_1 \in \GL_{h_1}(\D)$ and $Y_2 \in \GL_{h_2}(\D)$ such that
  \begin{align*}
    Y_1 B Y_2 =
    \begin{pmatrix}
      I_\ell & 0 \\
      0      & 0
    \end{pmatrix},
    \qquad
    Y_1 B' Y_2 =
    \begin{pmatrix}
      B'_{11} & B'_{12} \\
      B'_{21} & B'_{22}
    \end{pmatrix},
  \end{align*}
  where $B'_{11}$ is of size $\ell \times \ell$ and $B'_{22}$ is of size $(h_1 - \ell) \times (h_2 - \ell)$.
  The kernel of $Y_1 B Y_2$ consists of the vectors supported on the last $h_2 - \ell$ coordinates,
  and its image consists of the vectors supported on the first $\ell$ coordinates.
  Hence the desired inclusion is equivalent to $B'_{22} = 0$.

  Suppose, to the contrary, that $B'_{22} \neq 0$.
  Since $\emb$ is $\K$-linear and multiplicative, for $\lambda \in \K$ we have
  \begin{align*}
    \emb(Y_1)\, \emb(B + \lambda B')\, \emb(Y_2) =
    \begin{pmatrix}
      I_{d\ell} + \lambda \emb(B'_{11}) & \lambda \emb(B'_{12}) \\
      \lambda \emb(B'_{21})       & \lambda \emb(B'_{22})
    \end{pmatrix}.
  \end{align*}
  Since $\emb$ is injective, we have $\emb(B'_{22}) \neq 0$; we fix a nonzero entry of $\emb(B'_{22})$.
  Consider the minor of the matrix above with respect to the first $d\ell$ rows and columns together with the row and the column of this entry.
  It is of the form
  \begin{align*}
    f(\lambda) \coloneqq \det
    \begin{pmatrix}
      I_{d\ell} + \lambda \emb(B'_{11}) & \lambda \ast   \\
      \lambda \ast    & \lambda \gamma
    \end{pmatrix},
  \end{align*}
  where $\gamma \in \E$ denotes the fixed nonzero entry,
  and the two symbols $\ast$ denote some column vector in $\E^{d\ell}$ and some row vector in $\E^{1 \times d\ell}$, respectively.
  We regard $f$ as a polynomial in $\lambda$ over $\E$.
  Since every entry of the last row of the matrix above is a multiple of $\lambda$, we have $f(\lambda) = \lambda \tilde{f}(\lambda)$,
  where $\tilde{f}(\lambda)$ is the determinant of the matrix obtained by replacing the last row $(\lambda \ast \ \ \lambda \gamma)$ with $(\ast \ \ \gamma)$.
  Substituting $\lambda = 0$, we obtain $\tilde{f}(0) = \gamma$,
  since the resulting matrix is block lower triangular with diagonal blocks $I_{d\ell}$ and $\gamma$.
  Hence the coefficient of $\lambda$ in $f$ is $\gamma \neq 0$,
  and $f$ is a nonzero polynomial, which has finitely many roots in $\E$.
  Since $\K \subseteq \E$ is infinite, there is $\lambda_0 \in \K$ with $f(\lambda_0) \neq 0$.
  Then, by \Cref{lem:splitting-rank}, we obtain
  \begin{align*}
    d \rank_\D (B + \lambda_0 B') = \rank_\E \emb(B + \lambda_0 B') \geq d\ell + 1 > d\ell = d \rank_\D B,
  \end{align*}
  that is, $\rank_\D (B + \lambda_0 B') > \rank_\D B$, which contradicts the assumption.
\end{proof}

\section{Proof of the main theorem and consequences}\label{sec:proof}
This section is devoted to the proof of \Cref{thm:main}.
Let $\F$, $r$, and $A = A_0 x_0 + \sum_{i=1}^m A_i x_i$ be as in \Cref{thm:main}.
We fix an arbitrary integer $d \geq r$ and prove \Cref{thm:main} for $k = d$.

\subsection{Setup}\label{subsec:setup}
By \Cref{prop:division-algebra}, we take field extensions $\F \subseteq \K \subseteq \E$, a skew field $\D$, and an isomorphism $\spl$ such that
\begin{align*}
  \abs{\K} = \infty, \qquad Z(\D) = \K, \qquad N \coloneqq \dim_\K \D = d^2, \qquad \spl \colon \D \otimes_\K \E \to \Mat_d(\E).
\end{align*}
Let $\emb \colon \D \to \Mat_d(\E)$ be the embedding defined in \Cref{subsec:division-algebra}.
Since $d \geq r$, we have
\begin{align}\label{eq:numerical-condition}
  N = d^2 > r(r-1).
\end{align}
We regard $A_0, A_1, \dots, A_m$ also as matrices over $\K$ and over $\D$.
The rank of each $A_i$ does not depend on whether it is taken over $\F$, $\K$, or $\D$ (\Cref{lem:scalar-extension}~(3));
in particular, $\rank_\K A_i \leq r$ for $i \in [m]$.

For $W \leq \D^b$, we define
\begin{align}\label{eq:calA}
  \calA(W) \coloneqq \sum_{i = 1}^m A_i \D W = \spn_\D \Set{A_i(pw)}{i \in [m],\ p \in \D,\ w \in W} \leq \D^a.
\end{align}
Note that $A_0$ is not involved in the definition of $\calA$.
Moreover, $\calA(W)$ is two-sided:
Since $A_i$ has entries in the center, we have
$t \rbra{A_i(pw)p'} = A_i((t p)w)p' \in \calA(W)$
for all $t, p, p' \in \D$ and $w \in W$.

For a tuple $q = (q_1, \dots, q_m) \in \D^m$, we consider the matrix
\begin{align*}
  A(q) \coloneqq A_0 + \sum_{i=1}^m A_i q_i \in \D^{a \times b}
\end{align*}
together with the subspaces $V$ and $U$ and the integers $\rho$ and $\tau$ defined by
\begin{align}\label{eq:tuple-objects}
  \begin{aligned}
    V    & \coloneqq \calA(\ker_\D A(q)) \leq \D^a, & \qquad U    & \coloneqq A(q)^{-1}(V) = \Set{u \in \D^b}{A(q)u \in V} \leq \D^b, \\
    \rho & \coloneqq \rank_\D A(q),               & \qquad \tau & \coloneqq \dim_\D V.
  \end{aligned}
\end{align}
Here, $q_1, \dots, q_m$ are not formal variables but actual elements of the skew field $\D$;
that is, $A(q)$ is the matrix over $\D$ obtained from $A$ by substituting $1$ for $x_0$ and $q_i$ for $x_i$ for each $i \in [m]$.
A tuple $q \in \D^m$ is said to be \emph{lexicographically maximal}
if $\rho$ is maximum among all tuples in $\D^m$, and, subject to this, $\tau$ is maximum.
Since $\rho$ and $\tau$ are bounded nonnegative integers, a lexicographically maximal tuple exists.

We have $\ker_\D A(q) \leq U$, and hence $V = \calA(\ker_\D A(q)) \leq \calA(U)$.
Since $V$ is two-sided, \Cref{lem:central-descent} implies
\begin{align}\label{eq:W-descent}
  V = V_0 \otimes_\K \D, \qquad \dim_\K V_0 = \tau,
\end{align}
where $V_0 \coloneqq V \cap \K^a$.
The objects $V$, $U$, $V_0$, $\rho$, and $\tau$ depend on the tuple $q$, although we suppress it in the notation.

The proof proceeds as follows.
In \Cref{subsec:reduction}, we show that \Cref{thm:main} follows if some tuple $q$ satisfies the inclusion $V \leq \im_\D A(q)$ and the \emph{stability} $\calA(U) = V$ (\Cref{prop:reduction}),
and we state that every lexicographically maximal tuple satisfies both of them (\Cref{prop:lex-max}).
The subsequent two subsections are devoted to the proof of \Cref{prop:lex-max}.
In \Cref{subsec:maximality}, we prove the inclusion $V \leq \im_\D A(q)$.
In \Cref{subsec:stability}, we rephrase the stability $\calA(U) = V$ as the vanishing of certain vectors, and then prove it.
Finally, in \Cref{subsec:general-bound}, from \Cref{thm:main} we derive \Cref{cor:matrix-space} and a known bound for general linear symbolic matrices.

\subsection{Reduction to two conditions on a tuple}\label{subsec:reduction}
We first establish a lower bound on the rank of the $d$-th blow-up in terms of $\rho$, which holds for an arbitrary tuple.

\begin{lemma}\label{lem:lower-bound}
  Let $q \in \D^m$ be an arbitrary tuple, and let $\rho = \rank_\D A(q)$ as in~\eqref{eq:tuple-objects}.
  Then $\rank_{\F(\bX^{(d)})} A^{(d)} \geq d\rho$.
\end{lemma}
\begin{proof}
  Since the entries of $A_0, A_1, \dots, A_m$ belong to $\F \subseteq \K$ and $\emb(\kappa) = \kappa I_d$ for $\kappa \in \K$, we have
  \begin{align*}
    \emb(A(q)) = A_0 \otimes I_d + \sum_{i=1}^m A_i \otimes \emb(q_i) \in \E^{da \times db},
  \end{align*}
  that is, $\emb(A(q))$ is obtained from $A^{(d)}$ by substituting $I_d$ for $X_0^{(d)}$ and $\emb(q_i) \in \Mat_d(\E)$ for $X_i^{(d)}$ for each $i \in [m]$.
  By \Cref{lem:splitting-rank}, we have $\rank_\E \emb(A(q)) = d\rho$.
  Hence $\emb(A(q))$ has a nonzero minor of size $d\rho$.
  Consider the minor $f \in \F[\bX^{(d)}]$ of $A^{(d)}$ with the same rows and columns.
  The substitution above defines a ring homomorphism $\F[\bX^{(d)}] \to \E$, which sends $f$ to this nonzero minor of $\emb(A(q))$.
  Hence $f$ is a nonzero polynomial,
  and $A^{(d)}$ has a nonzero minor of size $d\rho$ over $\F(\bX^{(d)})$.
\end{proof}

The next proposition states that, if the tuple satisfies the following two conditions,
then the lower bound in \Cref{lem:lower-bound} is tight and $\rho$ coincides with the nc-rank of $A$.

\begin{proposition}\label{prop:reduction}
  Let $q \in \D^m$ be an arbitrary tuple, and let $V$, $U$, and $\rho$ be as in~\eqref{eq:tuple-objects}.
  Suppose that
  \begin{align}\label{eq:two-conditions}
    V \leq \im_\D A(q) \qquad \text{and} \qquad \calA(U) = V.
  \end{align}
  Then $\ncrk_\F A = \rho$, and~\eqref{eq:main} holds for $k = d$.
\end{proposition}
\begin{proof}
  It suffices to show $\ncrk_\K A \leq \rho$.
  Indeed, since the rank of a matrix does not change under a field extension (\Cref{lem:scalar-extension}~(3)),
  we have $\rank_{\F(\bX^{(k)})} A^{(k)} = \rank_{\K(\bX^{(k)})} A^{(k)}$ for every $k \geq 1$,
  and hence $\ncrk_\F A = \ncrk_\K A$ by \Cref{thm:blow-up}.
  Therefore, by \Cref{lem:lower-bound} and \Cref{thm:blow-up}, we obtain
  \begin{align*}
    d\rho \leq \rank_{\F(\bX^{(d)})} A^{(d)} \leq d \ncrk_\F A = d \ncrk_\K A \leq d\rho,
  \end{align*}
  which implies $\ncrk_\F A = \rho$ and~\eqref{eq:main} for $k = d$.

  To show $\ncrk_\K A \leq \rho$, we apply \Cref{thm:fortin-reutenauer} to $A$ regarded as a linear symbolic matrix over $\K$
  with the pair $(U_0, V_0)$, where $U_0 \coloneqq U \cap \K^b$ and $V_0 = V \cap \K^a$.
  For this, we show that $A_i U_0 \leq V_0$ holds for all $i \in [0, m]$ and that $\dim_\K U_0 - \dim_\K V_0 = b - \rho$;
  then \Cref{thm:fortin-reutenauer} yields $\ncrk_\K A \leq b - \dim_\K U_0 + \dim_\K V_0 = \rho$.

  We first show $A_i U \leq V$ for all $i \in [0, m]$.
  Let $u \in U$.
  By~\eqref{eq:calA} and the assumption $\calA(U) = V$, we have $A_i u \in V$ and $A_i(q_i u) \in V$ for all $i \in [m]$.
  Since $A(q)u \in V$ by the definition~\eqref{eq:tuple-objects} of $U$, we obtain $A_0 u = A(q)u - \sum_{i=1}^m A_i(q_i u) \in V$.
  Thus, $A_i U \leq V$ holds for all $i \in [0, m]$.
  Since $A_i$ has entries in $\K$, this implies $A_i U_0 \leq V \cap \K^a = V_0$ for all $i \in [0, m]$.

  We next show that $U$ is two-sided.
  Let $u \in U$ and $t \in \D$.
  Since $A_0$ has entries in the center, we have $A_0(t u) = t (A_0 u)$, and hence $A(q)(t u) = t (A_0 u) + \sum_{i=1}^m A_i((q_i t) u)$.
  The first term belongs to $V$, since $A_0 u \in V$ as shown above and $V$ is two-sided.
  Each $A_i((q_i t) u)$ belongs to $\calA(U) = V$ by~\eqref{eq:calA}.
  Hence $A(q)(t u) \in V$, that is, $t u \in U$.
  Thus $U$ is two-sided, and we have $U = U_0 \otimes_\K \D$ by \Cref{lem:central-descent}.
  Recall also from~\eqref{eq:W-descent} that $V = V_0 \otimes_\K \D$.

  Finally, by the assumption $V \leq \im_\D A(q)$, the restriction $A(q)|_U \colon U \to V$ is surjective, and its kernel is $\ker_\D A(q)$.
  Hence
  \begin{align*}
    \dim_\K U_0 - \dim_\K V_0 = \dim_\D U - \dim_\D V = \dim_\D \ker_\D A(q) = b - \rho,
  \end{align*}
  where the first equality follows from \Cref{lem:central-descent} and the last equality from the rank-nullity theorem over $\D$.
  This completes the proof.
\end{proof}

The two conditions~\eqref{eq:two-conditions} in \Cref{prop:reduction} are satisfied by every lexicographically maximal tuple.

\begin{proposition}\label{prop:lex-max}
  Let $q \in \D^m$ be a lexicographically maximal tuple, and let $V$ and $U$ be the subspaces defined in~\eqref{eq:tuple-objects}.
  Then the following hold.
  \begin{enumerate}[label=\textup{(\arabic*)}]
    \item $V \leq \im_\D A(q)$.
    \item $\calA(U) = V$.
  \end{enumerate}
\end{proposition}

We prove \Cref{prop:lex-max} in \Cref{subsec:maximality,subsec:stability}.
Before that, we complete the proof of \Cref{thm:main} by assuming \Cref{prop:lex-max}.

\begin{proof}[Proof of \Cref{thm:main}]
  Take a lexicographically maximal tuple $q \in \D^m$.
  By \Cref{prop:lex-max}, the two conditions in~\eqref{eq:two-conditions} hold for this tuple.
  Hence \Cref{thm:main} for $k = d$ follows from \Cref{prop:reduction}.
  Since $d \geq r$ is arbitrary, the theorem holds for all $k \geq r$.
\end{proof}

\begin{remark}\label{rem:wong-sequence}
  The subspaces $V$ and $U$ and the two conditions in~\eqref{eq:two-conditions} are analogous to the second Wong sequence of Ivanyos, Karpinski, Qiao, and Santha~\cite{Ivanyos2015}.
  For a matrix $B$ in a matrix space $\calM$ over a field,
  their \emph{second Wong sequence} is the sequence of subspaces defined by $W_0 \coloneqq \cbra{0}$ and $W_{i+1} \coloneqq \calM(B^{-1}(W_i))$ for $i \geq 0$,
  where $B^{-1}(W)$ denotes the preimage of $W$ under $B$,
  and $\calM(W)$ the subspace spanned by the images of $W$ under the matrices in $\calM$.
  The limit of this sequence is contained in the image of $B$ if and only if
  the rank of $B$ is certified to be equal to the nc-rank by a pair of subspaces as in \Cref{thm:fortin-reutenauer}~\cite[Lemma~9]{Ivanyos2015}.
  With $A(q)$ and $\calA$ in place of $B$ and $\calM$, the first two members of the sequence are $V = \calA(\ker_\D A(q))$ and $\calA(U)$,
  and the two conditions state that $V$ is contained in the image of $A(q)$ and that the sequence stabilizes at $V$.
\end{remark}

\subsection{The inclusion for a lexicographically maximal tuple}\label{subsec:maximality}
In this subsection, we derive consequences of the lexicographic maximality and prove \Cref{prop:lex-max}~(1).
To apply \Cref{lem:max-rank} to the lexicographic maximality, we combine the two quantities $\rho$ and $\tau$ into the rank of a single matrix as follows.

We fix an integer $c$ such that $c > a$ and $c > mN$, say
$c \coloneqq mN + a + 1$.
We fix a $\K$-basis $e_1, e_2, \dots, e_N$ of $\D$, and define matrices $L_1, \dots, L_c \in \D^{a \times b}$ as follows:
$L_1, \dots, L_{mN}$ are the matrices $A_j e_\alpha$ for $(j, \alpha) \in [m] \times [N]$ arranged in a fixed order, and $L_k \coloneqq 0$ for $mN < k \leq c$.
These matrices do not depend on the tuple $q$.
For every $W \leq \D^b$, we have
\begin{align}\label{eq:generators-L}
  \calA(W) = \sum_{j=1}^m \sum_{\alpha=1}^N A_j e_\alpha W = \sum_{k=1}^c L_k W,
\end{align}
since every $p \in \D$ admits a $\K$-linear combination $p = \sum_{\alpha} \kappa_\alpha e_\alpha$,
and $A_j(pw) = \sum_{\alpha} A_j e_\alpha (w \kappa_\alpha)$ with $w\kappa_\alpha \in W$ for all $j \in [m]$ and $w \in W$.

For a matrix $M \in \D^{a \times b}$, we define the block matrix
\begin{align*}
  P(M) \coloneqq
  \begin{pmatrix}
    M     & 0      & \cdots & 0      \\
    0      & M     & \cdots & 0      \\
    \vdots & \vdots & \ddots & \vdots \\
    0      & 0      & \cdots & M     \\
    L_1    & L_2    & \cdots & L_c
  \end{pmatrix}
  \in \D^{(ca + a) \times cb},
\end{align*}
where the diagonal part consists of $c$ copies of $M$.
In the following, we denote by $(w_1, w_2, \dots, w_c)$ the column vector obtained by concatenating column vectors $w_1, w_2, \dots, w_c$ vertically in this order,
and we refer to $w_k$ as its \emph{$k$-th block}.
For $w_1, w_2, \dots, w_c \in \D^b$, the matrix $P(M)$ maps the vector $(w_1, w_2, \dots, w_c) \in \D^{cb}$ to the vector $(M w_1, M w_2, \dots, M w_c, \sum_{k=1}^c L_k w_k) \in \D^{ca + a}$.

\begin{lemma}\label{lem:augmented-rank}
  Let $q \in \D^m$ be a lexicographically maximal tuple.
  Then, for every $q' \in \D^m$, we have
  \begin{align*}
    \rank_\D P(A(q')) \leq \rank_\D P(A(q)).
  \end{align*}
\end{lemma}
\begin{proof}
  We first show that, for every $M \in \D^{a \times b}$,
  \begin{align}\label{eq:rank-formula}
    \rank_\D P(M) = c \rank_\D M + \dim_\D \calA(\ker_\D M).
  \end{align}
  Consider the projection of $\im_\D P(M)$ to the first $c$ blocks, i.e., the map
  \begin{align*}
    \im_\D P(M) \ni \rbra*{Mw_1, \dots, Mw_c, \sum_{k=1}^c L_k w_k} \mapsto (Mw_1, \dots, Mw_c) \in \D^{ca},
  \end{align*}
  where $w_1, \dots, w_c \in \D^b$.
  Its image is $(\im_\D M)^c$, whose dimension is $c \rank_\D M$.
  Its kernel is identified with the set of the vectors $\sum_{k=1}^c L_k w_k$ with $w_1, \dots, w_c \in \ker_\D M$,
  which is equal to $\sum_{k=1}^c L_k \ker_\D M = \calA(\ker_\D M)$ by~\eqref{eq:generators-L}.
  Thus, we obtain~\eqref{eq:rank-formula} by the rank-nullity theorem over $\D$.

  Let $\rho' \coloneqq \rank_\D A(q')$ and $\tau' \coloneqq \dim_\D \calA(\ker_\D A(q'))$.
  By~\eqref{eq:rank-formula}, we have $\rank_\D P(A(q')) = c\rho' + \tau'$ and $\rank_\D P(A(q)) = c\rho + \tau$.
  By the lexicographic maximality, we have $\rho' \leq \rho$, and $\tau' \leq \tau$ if $\rho' = \rho$.
  If $\rho' < \rho$, then, by $\tau' \leq a < c$, we have $c\rho' + \tau' \leq c(\rho - 1) + a < c\rho \leq c\rho + \tau$.
  If $\rho' = \rho$, then $c\rho' + \tau' = c\rho + \tau' \leq c\rho + \tau$.
  Hence $\rank_\D P(A(q')) \leq \rank_\D P(A(q))$ holds in both cases.
\end{proof}

For a matrix $B \in \D^{a \times b}$, let $\Delta_B \in \D^{(ca + a) \times cb}$ be the block matrix obtained from $P(B)$ by replacing the last block row $(L_1 \ L_2 \ \cdots \ L_c)$ with zero, i.e.,
\begin{align*}
  \Delta_B \coloneqq
  \begin{pmatrix}
    B     & 0      & \cdots & 0      \\
    0      & B     & \cdots & 0      \\
    \vdots & \vdots & \ddots & \vdots \\
    0      & 0      & \cdots & B      \\
    0      & 0      & \cdots & 0
  \end{pmatrix}.
\end{align*}
Then $P(M) + \lambda \Delta_B = P(M + \lambda B)$ for all $\lambda \in \K$.

\begin{lemma}\label{lem:kernel-image}
  Let $q \in \D^m$ be a lexicographically maximal tuple.
  Then, for every $i \in [m]$ and $p \in \D$, we have
  $\Delta_{A_i p} \ker_\D P(A(q)) \leq \im_\D P(A(q))$.
\end{lemma}
\begin{proof}
  For every $\lambda \in \K$, we have $A(q) + \lambda A_i p = A(q')$,
  where $q' \in \D^m$ is the tuple obtained from $q$ by replacing $q_i$ with $q_i + \lambda p$.
  Hence, by \Cref{lem:augmented-rank}, we have
  $\rank_\D (P(A(q)) + \lambda \Delta_{A_i p}) = \rank_\D P(A(q) + \lambda A_i p) \leq \rank_\D P(A(q))$.
  Thus, the assertion follows from \Cref{lem:max-rank}.
\end{proof}

We are ready to prove $V \leq \im_\D A(q)$ (\Cref{prop:lex-max}~(1)).

\begin{proof}[Proof of \Cref{prop:lex-max}~\textup{(1)}]
  Let $q \in \D^m$ be a lexicographically maximal tuple.
  Let $i \in [m]$, $p \in \D$, and $u \in \ker_\D A(q)$.
  Take an index $k_0$ with $mN < k_0 \leq c$, so that $L_{k_0} = 0$.
  Let $w \in \D^{cb}$ be the vector whose $k_0$-th block is $u$ and whose other blocks are zero.
  Then $w \in \ker_\D P(A(q))$.
  By \Cref{lem:kernel-image}, there is a vector $w' \in \D^{cb}$ with $P(A(q)) w' = \Delta_{A_i p} w$.
  Comparing the $k_0$-th blocks, we obtain $A(q) w'_{k_0} = A_i(pu)$, where $w'_{k_0}$ denotes the $k_0$-th block of $w'$.
  Hence $A_i(pu) \in \im_\D A(q)$.
  Since $V = \calA(\ker_\D A(q))$ is spanned by such vectors by~\eqref{eq:calA}, we obtain $V \leq \im_\D A(q)$.
\end{proof}

\subsection{The stability for a lexicographically maximal tuple}\label{subsec:stability}
In this subsection, we prove \Cref{prop:lex-max}~(2).
We first rephrase the stability as the vanishing of certain vectors;
this rephrasing is valid for an arbitrary tuple $q$ satisfying the inclusion $V \leq \im_\D A(q)$.

Let $q \in \D^m$ be an arbitrary tuple.
We introduce a matrix that tests the membership in $V$.
Recall from~\eqref{eq:W-descent} that $V_0 = V \cap \K^a$ and $\dim_\K V_0 = \tau$.
Let $Q \in \K^{(a - \tau) \times a}$ be a matrix over $\K$ such that $\ker_\K Q = V_0$;
for instance, take a matrix representing the quotient map $\K^a \to \K^a / V_0 \simeq \K^{a - \tau}$.
By \Cref{lem:scalar-extension}~(3) and~\eqref{eq:W-descent}, we have
\begin{align}\label{eq:kernel-Q}
  \ker_\D Q = V.
\end{align}
That is, a vector $w \in \D^a$ belongs to $V$ if and only if $Q w = 0$.

Suppose that $V \leq \im_\D A(q)$.
Then there is a right $\D$-linear map $H \colon \D^a \to \D^b$ such that
\begin{align}\label{eq:section-H}
  A(q)Hv = v \qquad (v \in V).
\end{align}
Indeed, take a right $\D$-basis of $V$, choose a preimage under $A(q)$ for each basis vector,
extend the basis of $V$ to that of $\D^a$, and send the additional basis vectors to zero.
We identify $H$ with a matrix in $\D^{b \times a}$.
Note that $H$ does not preserve the left multiplication in general.

The following lemma rephrases the stability $\calA(U) = V$.

\begin{lemma}\label{lem:stability-criterion}
  Let $q \in \D^m$ be an arbitrary tuple, and let $V$ and $U$ be the subspaces defined in~\eqref{eq:tuple-objects}.
  Suppose that $V \leq \im_\D A(q)$,
  and let $Q$ be a matrix over $\K$ satisfying~\eqref{eq:kernel-Q} and $H$ a right $\D$-linear map satisfying~\eqref{eq:section-H}.
  Then $\calA(U) = V$ if and only if
  \begin{align}\label{eq:vanishing}
    Q A_j H A_i(pu) = 0 \qquad (i, j \in [m],\ p \in \D,\ u \in \ker_\D A(q)).
  \end{align}
\end{lemma}
\begin{proof}
  We first show $U = \ker_\D A(q) + HV$.
  For $u' \in \ker_\D A(q)$ and $v \in V$, we have $A(q)(u' + Hv) = v \in V$ by~\eqref{eq:section-H}, and hence $u' + Hv \in U$.
  Conversely, for $u \in U$, let $v \coloneqq A(q)u \in V$;
  then $A(q)(u - Hv) = v - v = 0$ by~\eqref{eq:section-H}, and hence $u = (u - Hv) + Hv \in \ker_\D A(q) + HV$.

  Since $\calA(W_1 + W_2) = \calA(W_1) + \calA(W_2)$ for $W_1, W_2 \leq \D^b$ by~\eqref{eq:calA} and $\calA(\ker_\D A(q)) = V$,
  we obtain $\calA(U) = V + \calA(HV)$.
  Thus $\calA(U) = V$ is equivalent to $\calA(HV) \leq V$.
  By~\eqref{eq:calA} and~\eqref{eq:kernel-Q}, this is equivalent to $Q A_j(t H v) = 0$ for all $j \in [m]$, $t \in \D$, and $v \in V$.
  Since $A_j$ and $Q$ have entries in the center, we have $Q A_j(t H v) = t\, Q A_j H v$;
  hence the condition is equivalent to $Q A_j H v = 0$ for all $j \in [m]$ and $v \in V$.
  Since $Q A_j H$ is right $\D$-linear and $V$ is spanned by the vectors $A_i(pu)$ with $i \in [m]$, $p \in \D$, and $u \in \ker_\D A(q)$ by~\eqref{eq:calA},
  this is further equivalent to~\eqref{eq:vanishing}.
\end{proof}

For a lexicographically maximal tuple $q \in \D^m$, we have $V \leq \im_\D A(q)$ by \Cref{prop:lex-max}~(1), which is already proved in \Cref{subsec:maximality};
hence a right $\D$-linear map $H$ satisfying~\eqref{eq:section-H} exists.
In the following, we also use the matrices $L_1, \dots, L_c$ and $P(A(q))$ defined in \Cref{subsec:maximality}, where $c$ is the number of blocks fixed there.

The following identity is the key consequence of the lexicographic maximality:
The vector $Q A_j H A_i(pu)$ with $u \in \ker_\D A(q)$ is obtained by multiplying $A_j u$ by a matrix over the center
that depends on neither $j$ nor $u$.

\begin{lemma}\label{lem:central-identity}
  Let $q \in \D^m$ be a lexicographically maximal tuple,
  $Q$ a matrix over $\K$ satisfying~\eqref{eq:kernel-Q},
  and $H$ a right $\D$-linear map satisfying~\eqref{eq:section-H}.
  Then, for each $i \in [m]$ and $p \in \D$, there is a matrix $Q_{i,p} \in \K^{(a - \tau) \times a}$ such that
  \begin{align*}
    Q_{i,p} A_j u = Q A_j H A_i (pu) \qquad (j \in [m],\ u \in \ker_\D A(q)).
  \end{align*}
\end{lemma}
\begin{proof}
  Fix $i \in [m]$ and $p \in \D$, and let $B \coloneqq A_i p \in \D^{a \times b}$.
  Note that $B u = A_i(pu) \in V$ for all $u \in \ker_\D A(q)$ by~\eqref{eq:calA}.

  Recall from~\eqref{eq:generators-L} that $V = \calA(\ker_\D A(q)) = \sum_{k=1}^c L_k \ker_\D A(q)$.
  We define a map $g \colon V \to \D^{a - \tau}$ by
  \begin{align*}
    g\rbra*{\sum_{k=1}^c L_k u_k} \coloneqq \sum_{k=1}^c Q L_k H B u_k \qquad (u_1, \dots, u_c \in \ker_\D A(q)).
  \end{align*}
  It suffices to show that $g$ is well-defined, is right $\D$-linear, and preserves the left multiplication, and that
  \begin{align}\label{eq:phi-generator}
    g(A_j(p'u)) = Q A_j (p' H B u)
  \end{align}
  holds for all $j \in [m]$, $p' \in \D$, and $u \in \ker_\D A(q)$.
  Indeed, since $V$ is two-sided, \Cref{lem:central-map} then provides a matrix $Q_{i,p} \in \K^{(a - \tau) \times a}$ such that $g(v) = Q_{i,p} v$ for all $v \in V$,
  and~\eqref{eq:phi-generator} with $p' = 1$ yields
  \begin{align*}
    Q_{i,p} A_j u = g(A_j u) = Q A_j H B u = Q A_j H A_i(pu)
  \end{align*}
  for all $j \in [m]$ and $u \in \ker_\D A(q)$, as required.
  In the following, we show these properties one by one.

  We first show that $g$ is well-defined, that is,
  the right-hand side depends only on the vector $\sum_{k} L_k u_k \in V$ and not on the choice of $u_1, \dots, u_c \in \ker_\D A(q)$ representing it.
  Suppose that $\sum_{k} L_k u_k = \sum_{k} L_k \tilde{u}_k$ for $u_k, \tilde{u}_k \in \ker_\D A(q)$.
  Then $u_k - \tilde{u}_k \in \ker_\D A(q)$ and $\sum_{k} L_k (u_k - \tilde{u}_k) = 0$,
  and the difference of the two values of the right-hand side is $\sum_{k} Q L_k H B (u_k - \tilde{u}_k)$.
  Hence the well-definedness is reduced to the following:
  If $\sum_{k} L_k u_k = 0$ for $u_1, \dots, u_c \in \ker_\D A(q)$, then $\sum_{k} Q L_k H B u_k = 0$.

  Suppose that $\sum_{k} L_k u_k = 0$ for $u_1, \dots, u_c \in \ker_\D A(q)$.
  Then the vector $(u_1, \dots, u_c) \in \D^{cb}$ belongs to $\ker_\D P(A(q))$.
  By \Cref{lem:kernel-image}, there are $w_1, \dots, w_c \in \D^b$ with
  \begin{align}
    &A(q) w_k = B u_k \quad (k = 1, \dots, c),\label{eq:kernel-image:1}\\
    &\sum_{k=1}^c L_k w_k = 0.\label{eq:kernel-image:2}
  \end{align}
  Since $B u_k \in V$, we have $A(q)HBu_k = Bu_k$ by~\eqref{eq:section-H},
  implying that $0 = Bu_k - A(q)HBu_k = A(q)(w_k - HBu_k)$ by~\eqref{eq:kernel-image:1}.
  Hence, we obtain $u'_k \coloneqq w_k - HBu_k \in \ker_\D A(q)$.
  Therefore,
  \begin{align*}
    \sum_{k=1}^c L_k H B u_k = \sum_{k=1}^c L_k w_k - \sum_{k=1}^c L_k u'_k = - \sum_{k=1}^c L_k u'_k \in \sum_{k=1}^c L_k \ker_\D A(q) = V
  \end{align*}
  holds, where the second equality follows from~\eqref{eq:kernel-image:2}.
  Thus, $\sum_{k} Q L_k H B u_k = 0$ by~\eqref{eq:kernel-Q}, as required.

  The right $\D$-linearity of $g$ follows from that of $H$, $B$, $L_k$, and $Q$,
  since $u_k + \tilde{u}_k$ and $u_k t$ belong to $\ker_\D A(q)$ for all $u_k, \tilde{u}_k \in \ker_\D A(q)$ and $t \in \D$.

  We next show~\eqref{eq:phi-generator}.
  Let $j \in [m]$, $p' \in \D$, and $u \in \ker_\D A(q)$,
  and write $p' = \sum_{\alpha} \kappa_\alpha e_\alpha$ with $\kappa_\alpha \in \K$.
  Then $A_j(p'u) = \sum_{\alpha} (A_j e_\alpha)(u \kappa_\alpha)$ with $u\kappa_\alpha \in \ker_\D A(q)$.
  By the definition of $g$, the right linearity of $H$ and $B$, and the centrality of $\kappa_\alpha$, we have
  \begin{align*}
    g(A_j(p'u))
    = \sum_{\alpha} Q A_j e_\alpha H B (u \kappa_\alpha)
    = \sum_{\alpha} Q A_j e_\alpha (H B u) \kappa_\alpha
    = Q A_j (p' HBu).
  \end{align*}

  Finally, we show that $g$ preserves the left multiplication.
  For any $t \in \D$, since $A_j$ and $Q$ have entries in the center,~\eqref{eq:phi-generator} implies
  \begin{align*}
    g(t A_j(p'u)) = g(A_j((t p')u)) = Q A_j((t p') HBu) = t\, Q A_j(p'HBu) = t g(A_j(p'u)).
  \end{align*}
  Since $V$ is spanned as a right $\D$-vector space by the vectors $A_j(p'u)$ by~\eqref{eq:calA} and $g$ is right $\D$-linear,
  we obtain $g(t v) = t g(v)$ for all $t \in \D$ and $v \in V$.
\end{proof}

We are now ready to prove the stability $\calA(U) = V$.

\begin{proof}[Proof of \Cref{prop:lex-max}~\textup{(2)}]
  Let $q \in \D^m$ be a lexicographically maximal tuple.
  By \Cref{prop:lex-max}~(1), we have $V \leq \im_\D A(q)$;
  let $Q$ be a matrix over $\K$ satisfying~\eqref{eq:kernel-Q} and $H$ a right $\D$-linear map satisfying~\eqref{eq:section-H}.
  Fix $i, j \in [m]$ and $u \in \ker_\D A(q)$.
  By \Cref{lem:stability-criterion}, it suffices to show that $Q A_j H A_i(pu) = 0$ for all $p \in \D$.

  If $\im_\D A_j \leq V$,
  then $Q A_j = 0$ by~\eqref{eq:kernel-Q}, and the assertion is trivial.

  In the following, we assume $\im_\D A_j \not\leq V$.
  Let $r_i \coloneqq \rank_\K A_i$ and $r_j \coloneqq \rank_\K A_j$; both are at most $r$.
  Recall from~\eqref{eq:tuple-objects} that $\tau = \dim_\D V$.
  If $a - \tau = 0$, then $V = \D^a$, which contradicts the assumption; hence $a - \tau \geq 1$.

  Fix any $k \in [a - \tau]$, and define a map $G \colon \D \to \D$ by
  \begin{align}\label{eq:G}
    G(p) \coloneqq \rbra*{Q A_j H A_i (pu)}_k = \rbra*{Q_{i,p} A_j u}_k,
  \end{align}
  where $(w)_k$ denotes the $k$-th entry of a vector $w \in \D^{a - \tau}$
  and the last equality follows from \Cref{lem:central-identity}.
  The map $G$ is $\K$-linear, since so are the map $p \mapsto pu$ and every right $\D$-linear map in the composition $Q A_j H A_i$.

  Our goal is to show $G = 0$.
  To this end, we bound the dimension of the image of $G$ by using $A_j$ (\Cref{cl:image-bound}),
  and the number of terms of $G$ as a sum of two-sided multiplications by using $A_i$ (\Cref{cl:term-bound}),
  as follows.

  \begin{claim}\label{cl:image-bound}
    $\rank_\K G \leq r_j - 1$.
  \end{claim}
  \begin{proof}[Proof of \Cref{cl:image-bound}]
  Let $V_j \coloneqq \spn_\D \Set{A_j(pu')}{p \in \D,\ u' \in \ker_\D A(q)}$.
  As in the case of $V$, the subspace $V_j \leq \D^a$ is two-sided.
  Moreover, we have $V_j \leq V \cap \im_\D A_j \lneq \im_\D A_j$,
  where the last strict inclusion follows from the assumption $\im_\D A_j \not\leq V$.
  Since $\dim_\D \im_\D A_j = r_j$ by \Cref{lem:scalar-extension}~(3), \Cref{lem:central-descent} implies
  $V_j = (V_j \cap \K^a) \otimes_\K \D$ and $\tau_j \coloneqq \dim_\K (V_j \cap \K^a) = \dim_\D V_j \leq r_j - 1$.

  Take a $\K$-basis $\eta_1, \dots, \eta_{\tau_j}$ of $V_j \cap \K^a$;
  it is also a right $\D$-basis of $V_j$ by \Cref{lem:scalar-extension}~(1) and \Cref{lem:central-descent}.
  Since $A_j u \in V_j$, there are $t_1, \dots, t_{\tau_j} \in \D$ such that $A_j u = \sum_{\ell=1}^{\tau_j} \eta_\ell t_\ell$.
  For each $p \in \D$, by~\eqref{eq:G}, we have
  \begin{align*}
    G(p) = \rbra*{Q_{i,p} A_j u}_k = \sum_{\ell=1}^{\tau_j} \rbra*{Q_{i,p} \eta_\ell}_k t_\ell.
  \end{align*}
  Since $Q_{i,p}$ and $\eta_\ell$ have entries in $\K$, each coefficient $(Q_{i,p} \eta_\ell)_k$ belongs to $\K$.
  Hence the image of $G$ is contained in the $\K$-linear span of $t_1, \dots, t_{\tau_j}$, which is independent of $p$.
  Therefore, we obtain $\rank_\K G \leq \tau_j \leq r_j - 1$.
  \end{proof}

  \begin{claim}\label{cl:term-bound}
    The map $G$ can be expressed as a sum of at most $r_i$ two-sided multiplications, that is,
    $G(p) = \sum_{\ell=1}^{r_i} y_\ell p z_\ell$ holds for all $p \in \D$
    for some $y_1, \dots, y_{r_i}, z_1, \dots, z_{r_i} \in \D$.
  \end{claim}
  \begin{proof}[Proof of \Cref{cl:term-bound}]
  Take a rank decomposition $A_i = \sum_{\ell=1}^{r_i} \xi_\ell \zeta_\ell$ over $\K$,
  where $\xi_\ell \in \K^a$ is a column vector and $\zeta_\ell \in \K^{1 \times b}$ is a row vector.
  Let $z_\ell \coloneqq \zeta_\ell u \in \D$.
  Since $\zeta_\ell$ has entries in the center, we have $\zeta_\ell(pu) = p \zeta_\ell u = p z_\ell$.
  Thus, for every $p \in \D$, we obtain $A_i(pu) = \sum_{\ell=1}^{r_i} \xi_\ell p z_\ell$.
  By the right $\D$-linearity of $Q A_j H$ and~\eqref{eq:G},
  \begin{align*}
    G(p) = \rbra*{Q A_j H A_i (pu)}_k = \sum_{\ell=1}^{r_i} \rbra*{Q A_j H \xi_\ell}_k p z_\ell
  \end{align*}
  holds for $p \in \D$.
  Letting $y_\ell \coloneqq \rbra*{Q A_j H \xi_\ell}_k \in \D$ for $\ell \in [r_i]$, we obtain the required expression $G(p) = \sum_{\ell=1}^{r_i} y_\ell p z_\ell$.
  \end{proof}

  We now show $G = 0$.
  By~\eqref{eq:numerical-condition} and \Cref{cl:image-bound}, we have
  \begin{align*}
    r_i \rank_\K G \leq r_i (r_j - 1) \leq r(r-1) < N.
  \end{align*}
  Hence, \Cref{lem:operator-bound} and \Cref{cl:term-bound} imply $G = 0$.
\end{proof}

\subsection{Consequences}\label{subsec:general-bound}
We first derive \Cref{cor:matrix-space} from \Cref{thm:main}.
Recall that a \emph{matrix space} is an $\F$-subspace of $\F^{a \times b}$.
For a matrix space $\calM \leq \F^{a \times b}$ and a positive integer $k$, its \emph{$k$-th tensor blow-up} is the matrix space $\calM \otimes_\F \F^{k \times k} \leq \F^{ak \times bk}$~\cite[Definition~1.2]{Ivanyos2018} (see also \cite[Definition~1.5]{Derksen2018}),
that is, the subspace spanned by the matrices $B \otimes T$ with $B \in \calM$ and $T \in \F^{k \times k}$.
The $k$-th blow-up of a linear symbolic matrix whose coefficients span $\calM$ is a generic matrix in the $k$-th tensor blow-up of $\calM$.

The ranks of the blow-ups and the nc-rank depend only on the matrix space spanned by the coefficients.
This is a standard fact~\cite[Corollary~1 and Lemma~1]{Fortin2004} (see also~\cite[Lemma~1.4]{Derksen2018}),
which is used, e.g., by Oki and Soma~\cite[Section~6.3.5]{Oki2025} in their treatment of the hidden fractional linear matroid parity problem.
We state it for an arbitrary field in terms of the ranks over rational function fields,
and include a proof for completeness.

\begin{lemma}\label{lem:matrix-space}
  Let $A = \sum_{i=0}^m A_i x_i$ and $B = \sum_{j=0}^{m'} B_j x_j$ be linear symbolic matrices over a field $\F$ with $A_i, B_j \in \F^{a \times b}$.
  If $A_0, A_1, \dots, A_m$ and $B_0, B_1, \dots, B_{m'}$ span the same matrix space,
  then the ranks of $A^{(k)}$ and $B^{(k)}$ are the same for every $k \geq 1$,
  and $\ncrk_\F A = \ncrk_\F B$.
  Here, the rank of $B^{(k)} = \sum_{j=0}^{m'} B_j \otimes X_j^{(k)}$ is taken over the rational function field over $\F$ in the entries of $X_0^{(k)}, X_1^{(k)}, \dots, X_{m'}^{(k)}$.
\end{lemma}
\begin{proof}
  Let $\F'$ be an infinite extension field of $\F$, and let $k \geq 1$.
  Since the rank of a matrix does not change under a field extension,
  and since a polynomial over the infinite field $\F'$ is nonzero if and only if it does not vanish at some point,
  the rank of $A^{(k)}$ is equal to the maximum rank of $\sum_{i=0}^m A_i \otimes T_i$ over all $T_0, T_1, \dots, T_m \in (\F')^{k \times k}$;
  apply this to the minors of $A^{(k)}$.
  The set of the matrices of the form $\sum_{i=0}^m A_i \otimes T_i$ is the $k$-th tensor blow-up of the matrix space over $\F'$ spanned by $A_0, A_1, \dots, A_m$.
  Hence it depends only on this matrix space, and the ranks of $A^{(k)}$ and $B^{(k)}$ are the same.
  The equality $\ncrk_\F A = \ncrk_\F B$ follows from \Cref{thm:blow-up}.
\end{proof}

\begin{proof}[Proof of \Cref{cor:matrix-space}]
  Let $B \coloneqq \sum_{j=0}^{m'} B_j x_j$ and $k \geq r$.
  By \Cref{thm:main} applied to $B$, the $k$-th blow-up of $B$ attains $\ncrk B$.
  Since $A_0, A_1, \dots, A_m$ and $B_0, B_1, \dots, B_{m'}$ span the same matrix space,
  the ranks of $A^{(k)}$ and $B^{(k)}$ are the same and $\ncrk A = \ncrk B$ by \Cref{lem:matrix-space}.
  Hence the $k$-th blow-up of $A$ attains $\ncrk A$, that is, \eqref{eq:main} holds.
\end{proof}

As mentioned in \Cref{sec:intro}, \Cref{thm:main} also yields an alternative proof of the bound of Derksen and Makam in the following form,
which requires no assumption on the ranks of the coefficients.
\begin{corollary}[{cf.~\cite[Theorem~1.8]{Derksen2017} and \cite[Section~2.1]{Huszar2021}}]\label{cor:general-bound}
  Let $A$ be a linear symbolic matrix over a field $\F$ of the form~\eqref{eq:linear-matrix}.
  Then, for every positive integer $k \geq \ncrk A - 1$, the $k$-th blow-up of $A$ attains the nc-rank, that is,
  \begin{align*}
    \ncrk A = \frac{1}{k} \rank_{\F(\bX^{(k)})} A^{(k)}.
  \end{align*}
  In particular, this holds for every positive integer $k \geq \min \cbra{a, b} - 1$.
\end{corollary}
\begin{proof}
  For every $i \in [0, m]$ and $k \geq 1$, we have
  $\rank A_i \leq (\rank_{\F(\bX^{(k)})} A^{(k)})/k \leq \ncrk A$,
  where the first inequality is obtained by substituting $I_k$ for $X_i^{(k)}$ and the zero matrix for the others,
  and the second one follows from \Cref{thm:blow-up}.
  If $\rank A_i = \ncrk A$ for some $i$, then the $k$-th blow-up attains the nc-rank for every $k \geq 1$.
  Otherwise, we have $\rank A_i \leq \ncrk A - 1$ for all $i$.
  Then $\ncrk A \geq 2$, since otherwise all $A_i$ are zero and hence $\ncrk A = 0 = \rank A_0$, a contradiction.
  Hence the $k$-th blow-up attains the nc-rank for every $k \geq \ncrk A - 1$ by \Cref{thm:main} applied with $r = \ncrk A - 1$.
  The last assertion follows from $\ncrk A \leq \min \cbra{a, b}$.
\end{proof}

\section{Applications}\label{sec:application}
In this section, we give two applications of \Cref{thm:main}.
In \Cref{subsec:randomized}, we present a fast randomized algorithm for computing the nc-rank.
In \Cref{subsec:newton}, we show that the hierarchy of polytopes of Hirai~\cite{Hirai2026} stabilizes at the $r$-th level.

\subsection{Fast randomized computation of the noncommutative rank}\label{subsec:randomized}
Let $A$ be a linear symbolic matrix over a field $\F$ of the form~\eqref{eq:linear-matrix} with $A_0, A_1, \dots, A_m \in \F^{a \times b}$
that satisfies the assumption of \Cref{cor:matrix-space} for a positive integer $r$.
Then the $r$-th blow-up of $A$ attains the nc-rank, that is, $\ncrk_\F A = (\rank_{\F(\bX^{(r)})} A^{(r)})/r$.
The rank of $A^{(r)}$ is attained with a constant probability
by substituting random elements of a sufficiently large finite subset $S \subseteq \F$ for the variables.
This yields the simple randomized algorithm for computing the nc-rank described in \Cref{alg:randomized}:
It substitutes $r \times r$ matrices $T_0, T_1, \dots, T_m$ with random entries in $S$ for $X_0^{(r)}, X_1^{(r)}, \dots, X_m^{(r)}$ in $A^{(r)}$,
computes the rank of the resulting matrix $A(\bT) \coloneqq \sum_{i=0}^m A_i \otimes T_i \in \F^{ar \times br}$ over $\F$,
and outputs the rank divided by $r$ and rounded up.

\begin{algorithm}[t]
  \caption{Randomized algorithm for the nc-rank}\label{alg:randomized}
  \begin{algorithmic}[1]
    \Require A linear symbolic matrix $A = A_0 x_0 + A_1 x_1 + \dots + A_m x_m$ with $A_0, A_1, \dots, A_m \in \F^{a \times b}$,
    a positive integer $r$, and a finite subset $S \subseteq \F$ with $\abs{S} \geq 2 r \min \cbra{a, b}$.
    \Ensure An integer that is at most $\ncrk_\F A$.
    \State Choose all the entries of $r \times r$ matrices $T_0, T_1, \dots, T_m$ independently and uniformly at random from $S$.
    \State Construct the matrix $A(\bT) = \sum_{i=0}^m A_i \otimes T_i \in \F^{ar \times br}$.
    \State Compute $\rank_\F A(\bT)$, and output $\ceil{(\rank_\F A(\bT)) / r}$.
  \end{algorithmic}
\end{algorithm}

The following theorem shows the correctness and the running time of \Cref{alg:randomized}.

\begin{theorem}\label{thm:randomized}
  Let $A$ be a linear symbolic matrix over a field $\F$ that satisfies the assumption of \Cref{cor:matrix-space} for a positive integer $r$,
  and suppose that $\abs{\F} \geq 2 r \min \cbra{a, b}$.
  Then \Cref{alg:randomized}
  computes $\ncrk_\F A$ with probability at least $1/2$
  in $O(m a b r^2 + (nr)^\omega)$ arithmetic operations in $\F$,
  where $n \coloneqq \max \cbra{a, b}$ and $2 \leq \omega \leq 3$ is the matrix multiplication exponent.
\end{theorem}
\begin{proof}
  Let $\rho \coloneqq \ncrk_\F A$.
  By \Cref{cor:matrix-space} with $k = r$, we have $\rank_{\F(\bX^{(r)})} A^{(r)} = r \rho$.
  Since $A(\bT)$ is obtained from $A^{(r)}$ by a substitution, we have $\rank_\F A(\bT) \leq r \rho$,
  and hence the output is at most $\rho$ for every choice of $\bT$.
  Moreover, $A^{(r)}$ has a nonzero minor $f$ of size $r \rho$,
  which is a polynomial of total degree $r \rho$ in the entries of the matrices in $\bX^{(r)}$.
  If $f$ does not vanish at $\bT$, then $\rank_\F A(\bT) = r \rho$, and the output is $\rho$.
  By the Schwartz--Zippel lemma~\cite[Corollary~1]{Schwartz1980} (see also~\cite{Zippel1979}),
  a nonzero polynomial of total degree $\ell$ over a field vanishes with probability at most $\ell / \abs{S}$
  at a point whose coordinates are chosen independently and uniformly at random from a finite subset $S$ of the field.
  Hence $f$ vanishes at $\bT$ with probability at most $r \rho / \abs{S} \leq r \min \cbra{a, b} / \abs{S} \leq 1/2$.

  The matrix $A(\bT)$ is constructed with $O(m a b r^2)$ arithmetic operations,
  and its rank is computed with $O((nr)^\omega)$ arithmetic operations by regarding it as a square matrix of size $nr$ after appending zero rows or zero columns.
\end{proof}

If $r$ is a constant, then the number of arithmetic operations in \Cref{thm:randomized} is $O(m n^2 + n^\omega)$.
That is, the nc-rank is computed in time linear in the input size plus the time for computing the rank of a single $n \times n$ matrix over $\F$.

\paragraph{Bit complexity over the rationals.}
For the case $\F = \Q$, the bit lengths of the numbers appearing in \Cref{alg:randomized} are bounded by a polynomial in the input size;
a similar analysis for the fractional linear matroid parity problem is given in~\cite[Section~6.3.3]{Oki2025}.

For a matrix $B$ over $\Q$, let $\bitsize{B}$ denote its bit complexity,
that is, the sum of the bit lengths of the numerators and denominators of its entries,
and let $\bitsize{A} \coloneqq \max \Set{\bitsize{A_i}}{i \in [0, m]}$.

\begin{proposition}\label{prop:bit-complexity}
  Let $A$ be as in \Cref{thm:randomized} with $\F = \Q$.
  Then \Cref{alg:randomized} with $S = [0, 2 r \min \cbra{a, b} - 1]$
  computes $\ncrk_\Q A$ with probability at least $1/2$,
  and it can be implemented so that every number appearing in the computation is an integer of bit length $O(nr (\bitsize{A} + \log (m n r)))$.
\end{proposition}
\begin{proof}
  By multiplying each $A_i$ by the product of the denominators of its entries,
  we may assume that the entries of $A_0, A_1, \dots, A_m$ are integers of bit length at most $\bitsize{A}$.
  The output is $\ncrk_\Q A$ with probability at least $1/2$ by \Cref{thm:randomized}.
  Each entry of $A(\bT)$ is a sum of $m + 1$ products of an entry of $A_i$ and an entry of $T_i$,
  and hence is an integer of bit length $\beta = O(\bitsize{A} + \log (m n r))$.
  We compute the rank of $A(\bT)$ by the variant of Gaussian elimination due to Edmonds~\cite[Section~7]{Edmonds1967},
  in which every number appearing in the computation is a minor of $A(\bT)$~\cite[Theorem~2]{Edmonds1967}.
  A minor of size $k$ is a sum of $k!$ signed products of $k$ entries,
  and hence is an integer of bit length at most $k (\beta + \log k) + 1$.
  Since $k \leq nr$, this is $O(nr (\bitsize{A} + \log (m n r)))$.
\end{proof}

\paragraph{Extension to small finite fields.}
\Cref{alg:randomized} requires that $\F$ has at least $2 r \min \cbra{a, b}$ elements.
Even if $\F$ is a finite field with $\abs{\F} < 2 r \min \cbra{a, b}$,
the nc-rank can still be computed by a randomized algorithm:
We apply \Cref{alg:randomized} over an extension field of $\F$.
Take an extension field $\F'$ of $\F$ with $\abs{\F'} \geq 2 r \min \cbra{a, b}$.
The rank of a matrix over $\F(\bX^{(r)})$ does not change if it is regarded as a matrix over $\F'(\bX^{(r)})$,
by \Cref{lem:scalar-extension}~(3) applied to the field $\F'(\bX^{(r)})$ and its subfield $\F(\bX^{(r)})$;
hence we also have $\rank_{\F'(\bX^{(r)})} A^{(r)} = r \ncrk_\F A$.
Since the proof of \Cref{thm:randomized} uses only this equality and the Schwartz--Zippel lemma,
\Cref{alg:randomized} in which $S$ is taken from $\F'$ instead of $\F$ and the rank of $A(\bT)$ is computed over $\F'$
outputs $\ncrk_\F A$ with probability at least $1/2$.

Concretely, we can construct $\F'$ as follows.
Let $h$ be the smallest positive integer such that $\abs{\F}^h \geq 2 r \min \cbra{a, b}$.
We take a monic irreducible polynomial $g \in \F[z]$ of degree $h$,
and let $\F' \coloneqq \F[z] / (g)$ and $\theta \coloneqq z \bmod g$.
Such a polynomial $g$ can be found by the randomized algorithm of Shoup~\cite[Theorem~5.1]{Shoup1994},
whose expected number of arithmetic operations in $\F$ is bounded by a polynomial in $h$ and $\log \abs{\F}$.
Since $h = O(\log (nr))$ and $\abs{\F} < 2 r \min \cbra{a, b}$, this number is $(\log (nr))^{O(1)}$.
Then $\F'$ is a field with $\abs{\F}^h$ elements, and we take $S \coloneqq \F'$.
We represent each element of $\F'$ by its coefficient vector with respect to the basis $1, \theta, \dots, \theta^{h-1}$ over $\F$.
An addition in $\F'$ uses $O(h)$ arithmetic operations in $\F$.
A multiplication uses $O(h^2)$ arithmetic operations in $\F$, by multiplying the representing polynomials and reducing the product modulo $g$.
The inverse of a nonzero element uses $O(h^2)$ arithmetic operations in $\F$, by the extended Euclidean algorithm applied to its representing polynomial and $g$~\cite[Corollary~4.6]{Gathen2013}.
A uniformly random element of $\F'$ is obtained from $h$ independent uniformly random elements of $\F$ as its coefficients.
Therefore, apart from finding $g$, \Cref{alg:randomized} uses $h (m+1) r^2$ independent random elements of $\F$
and $O\rbra{h^2 \rbra{m a b r^2 + (nr)^\omega}}$ arithmetic operations in $\F$.
Since $h = O(\log (nr))$, this is $O\rbra{(m a b r^2 + (nr)^\omega) \log^2 (nr)}$,
which dominates the expected number of arithmetic operations for finding $g$.

\subsection{Hierarchy for the nc-independent set polytope}\label{subsec:newton}
Hirai~\cite{Hirai2026} introduced a hierarchy of polytopes associated with a linear symbolic matrix,
whose limit is a noncommutative counterpart of the Newton polytope of the minors.

Let $A$ be a linear symbolic matrix over a field $\F$ of the form~\eqref{eq:linear-matrix} with $A_0, A_1, \dots, A_m \in \F^{a \times b}$.
For $I \subseteq [a]$ and $J \subseteq [b]$ with $\abs{I} = \abs{J}$,
we denote by $A[I, J] \coloneqq \sum_{i=0}^m A_i[I, J] x_i$ the linear symbolic matrix obtained from $A$ by taking the rows in $I$ and the columns in $J$,
where $A_i[I, J]$ is the corresponding submatrix of $A_i$.
For a positive integer $k$, the $k$-th blow-up $A[I, J]^{(k)}$ is the submatrix of $A^{(k)}$ consisting of the blocks indexed by $I \times J$.
We regard an exponent vector of a monomial in the entries of the matrices in $\bX^{(k)}$ as a vector $\alpha = (\alpha_{i, j\ell})$ of nonnegative integers indexed by $(i, j, \ell) \in [0, m] \times [k] \times [k]$,
where $\alpha_{i, j\ell}$ is the exponent of the $(j, \ell)$-th entry of $X_i^{(k)}$.
Let $\proj_k(\alpha) \in \R^{m+1}$ be the vector whose coordinates are indexed by $[0, m]$ and given by
\begin{align*}
  \proj_k(\alpha)_i \coloneqq \frac{1}{k} \sum_{j, \ell \in [k]} \alpha_{i, j\ell} \qquad (i \in [0, m]).
\end{align*}
For a polynomial $f \in \F[\bX^{(k)}]$, we denote by $\supp f$ the set of the exponent vectors of the monomials appearing in $f$ with nonzero coefficients,
that is, $\supp f \coloneqq \Set*{\alpha}{f_\alpha \neq 0}$,
where $f_\alpha \in \F$ denotes the coefficient of the monomial $\prod_{i \in [0, m]} \prod_{j, \ell \in [k]} (X_i^{(k)})_{j\ell}^{\alpha_{i, j\ell}}$ in $f$.
Following Hirai~\cite[(2.7) and (3.1)]{Hirai2026}, we define
\begin{align*}
  \calP^{k}(A) &\coloneqq \conv \bigcup_{\substack{I \subseteq [a],\, J \subseteq [b] \\ \abs{I} = \abs{J}}} \proj_k\rbra*{\supp \det A[I, J]^{(k)}}, \\
  \calP^{\leq k}(A) &\coloneqq \conv \bigcup_{k'=1}^{k} \calP^{k'}(A), \qquad
  \calQ(A) \coloneqq \conv \bigcup_{k' \geq 1} \calP^{k'}(A),
\end{align*}
where the first union includes the pair $I = J = \emptyset$ (recall that the determinant of the empty matrix is $1$),
and $\conv$ denotes the convex hull in $\R^{m+1}$.
The polytope $\calQ(A)$ is called the \emph{nc-independent set polytope} of $A$~\cite[Section~3.2.4]{Hirai2025}.
Note that $\calP^{k}(A)$ and $\calP^{\leq k}(A)$ are polytopes, whereas $\calQ(A)$ is defined as the convex hull of an infinite union;
it is shown to be an integral polytope in~\cite[Theorem~3.17~(2)]{Hirai2025}.
The polytope $\calP^{1}(A)$, which is denoted by $\calP(A)$ and called the \emph{independent set polytope} of $A$ in~\cite[Section~3.2.4]{Hirai2025}, is the convex hull of the union of the Newton polytopes of all minors of $A$ regarded as a matrix in commuting variables $x_0, x_1, \dots, x_m$.

By definition, these polytopes form the hierarchy
$\calP^{1}(A) = \calP^{\leq 1}(A) \subseteq \calP^{\leq 2}(A) \subseteq \dots \subseteq \calQ(A)$.
The following theorem shows that this hierarchy stabilizes at the $r$-th level under the assumption of \Cref{thm:main}.

\begin{theorem}\label{thm:newton}
  Let $\F$, $r$, and $A$ be as in \Cref{thm:main}.
  Then, for every integer $k \geq r$, we have
  \begin{align*}
    \calP^{k}(A) = \calP^{\leq k}(A) = \calQ(A).
  \end{align*}
\end{theorem}

In particular, we have $\calQ(A) = \calP^{r}(A)$.
For $r = 1$, the nc-independent set polytope $\calQ(A)$ coincides with the independent set polytope $\calP^{1}(A)$.

To prove \Cref{thm:newton}, we use the results of~\cite{Hirai2025} on the degree of the Dieudonn\'e determinant of a linear symbolic matrix whose coefficients are rational functions in another variable $t$.
For a nonzero $f \in \F(t)$, written as $f = f_1 / f_2$ with polynomials $f_1, f_2 \in \F[t]$,
let $\deg f \coloneqq \deg f_1 - \deg f_2$, and let $\deg 0 \coloneqq -\infty$.
For a polynomial $f$ in commuting variables with coefficients in $\F(t)$,
let $\deg f$ denote the maximum of $\deg c$ over the coefficients $c$ of $f$.
That is, throughout this subsection, the degree is always taken with respect to $t$.
Let $\calO_\infty \coloneqq \Set{f \in \F(t)}{\deg f \leq 0}$, which is a subring of $\F(t)$.
For $f = f_1 / f_2 \in \calO_\infty$ with $\deg f_1 \leq \deg f_2$,
let $f(\infty) \in \F$ be the ratio of the coefficients of $t^{\deg f_2}$ in $f_1$ and $f_2$.
Then $f \mapsto f(\infty)$ is a ring homomorphism $\calO_\infty \to \F$.
For a matrix $R$ over $\calO_\infty$, we denote by $R(\infty)$ the matrix over $\F$ of the same size
obtained by applying this homomorphism to each entry of $R$.

Let $\ell$ be a positive integer,
and let $B = \sum_{i=0}^m B_i x_i$ be a linear symbolic matrix over $\F(t)$ with $B_0, B_1, \dots, B_m \in \F(t)^{\ell \times \ell}$.
We denote by $\deg \Det B \in \Z \cup \cbra{-\infty}$ the degree of the \emph{Dieudonn\'e determinant}~\cite{Dieudonne1943} of $B$,
where $B$ is regarded as a matrix over the skew field of rational functions in $t$ over the free skew field $\freefield{\F}$,
and $\deg \Det B \coloneqq -\infty$ if $B$ is singular over this skew field;
see~\cite[Section~2.2]{Hirai2019} and~\cite[Section~2.2]{Hirai2025} for the definition.
A pair $(Y_1, Y_2)$ of matrices in $\GL_\ell(\F(t))$ is said to be \emph{feasible} for $B$ if $Y_1 B_i Y_2 \in \calO_\infty^{\ell \times \ell}$ for all $i \in [0, m]$.
For a feasible pair $(Y_1, Y_2)$, we denote by $(Y_1 B Y_2)(\infty) \coloneqq \sum_{i=0}^m (Y_1 B_i Y_2)(\infty) x_i$ the resulting linear symbolic matrix over $\F$.

\begin{proposition}[{\cite[Theorem~2.4, Lemma~2.5, and Proposition~2.6]{Hirai2025}}]\label{prop:deg-det}
  Let $B$ be a linear symbolic matrix over $\F(t)$ with $\ell \times \ell$ coefficients as above.
  \begin{enumerate}[label=\textup{(\arabic*)}]
    \item $\deg \Det B$ is the infimum of $-\deg \det Y_1 - \deg \det Y_2$ over all feasible pairs $(Y_1, Y_2)$ for $B$.
    In addition, a feasible pair $(Y_1, Y_2)$ attains the infimum if and only if $\ncrk_\F (Y_1 B Y_2)(\infty) = \ell$.
    \item If a feasible pair $(Y_1, Y_2)$ for $B$ satisfies $\rank_{\ratfield{\F}} (Y_1 B Y_2)(\infty) = \ell$,
    then $\deg \det B = -\deg \det Y_1 - \deg \det Y_2$.
    \item $\deg \Det B = \max_{k \geq 1} (\deg \det B^{(k)})/k$.
  \end{enumerate}
\end{proposition}

The following proposition is a counterpart of \Cref{thm:main} for the degree of the Dieudonn\'e determinant:
The maximum in \Cref{prop:deg-det}~(3) is attained by every $k \geq r$ if all but one of the coefficients have rank at most $r$.

\begin{proposition}\label{prop:deg-det-attainment}
  Let $\ell$ and $r$ be positive integers,
  and let $B = \sum_{i=0}^m B_i x_i$ be a linear symbolic matrix over $\F(t)$ with $B_0, B_1, \dots, B_m \in \F(t)^{\ell \times \ell}$ such that $\rank_{\F(t)} B_i \leq r$ for all $i \in [m]$;
  the matrix $B_0$ is arbitrary.
  Then, for every integer $k \geq r$, the $k$-th blow-up of $B$ attains the degree of the Dieudonn\'e determinant, that is,
  \begin{align*}
    \deg \Det B = \frac{1}{k} \deg \det B^{(k)}.
  \end{align*}
\end{proposition}
\begin{proof}
  Let $k \geq r$ be an integer.
  If $\deg \Det B = -\infty$, then $\deg \det B^{(k)} = -\infty$ by \Cref{prop:deg-det}~(3).
  Suppose that $\deg \Det B > -\infty$.
  Since $-\deg \det Y_1 - \deg \det Y_2$ is an integer for every feasible pair,
  the infimum in \Cref{prop:deg-det}~(1) is attained by some feasible pair $(Y_1, Y_2)$ for $B$.
  For $i \in [0, m]$, let $C_i \coloneqq (Y_1 B_i Y_2)(\infty) \in \F^{\ell \times \ell}$,
  and let $C \coloneqq \sum_{i=0}^m C_i x_i = (Y_1 B Y_2)(\infty)$.
  Then $\ncrk_\F C = \ell$ by \Cref{prop:deg-det}~(1).

  We next show that $\rank_\F C_i \leq r$ for all $i \in [m]$.
  Let $i \in [m]$.
  Since $Y_1$ and $Y_2$ are nonsingular, we have $\rank_{\F(t)} Y_1 B_i Y_2 = \rank_{\F(t)} B_i \leq r$,
  and hence every minor of $Y_1 B_i Y_2$ of size larger than $r$ is zero.
  Since $f \mapsto f(\infty)$ is a ring homomorphism,
  each minor of $C_i$ is obtained by applying it to the corresponding minor of $Y_1 B_i Y_2$.
  Therefore, every minor of $C_i$ of size larger than $r$ is zero, and $\rank_\F C_i \leq r$.

  By \Cref{thm:main} applied to $C$, we have $\rank_{\F(\bX^{(k)})} C^{(k)} = k \ell$.
  We regard $B^{(k)}$ as a linear symbolic matrix over $\F(t)$ whose variables are the entries of the matrices in $\bX^{(k)}$;
  the coefficient of the $(j, j')$-th entry of $X_i^{(k)}$ is $B_i \otimes E_{jj'}$,
  where $E_{jj'}$ denotes the $k \times k$ matrix whose $(j, j')$-th entry is $1$ and whose other entries are $0$.
  Since $(Y_1 \otimes I_k)(B_i \otimes E_{jj'})(Y_2 \otimes I_k) = (Y_1 B_i Y_2) \otimes E_{jj'}$,
  the pair $(Y_1 \otimes I_k, Y_2 \otimes I_k)$ is feasible for $B^{(k)}$,
  and $\rbra{(Y_1 \otimes I_k) B^{(k)} (Y_2 \otimes I_k)}(\infty) = C^{(k)}$.
  Therefore, by \Cref{prop:deg-det}~(2) applied to $B^{(k)}$ and $\det (Y \otimes I_k) = (\det Y)^k$ for $Y \in \GL_\ell(\F(t))$, we obtain
  $\deg \det B^{(k)} = -k \deg \det Y_1 - k \deg \det Y_2 = k \deg \Det B$.
\end{proof}

We are now ready to prove \Cref{thm:newton}.

\begin{proof}[Proof of \Cref{thm:newton}]
  Since $\calP^{k}(A) \subseteq \calP^{\leq k}(A) \subseteq \calQ(A)$ holds for every positive integer $k$ by definition,
  it suffices to show $\calQ(A) \subseteq \calP^{k}(A)$ for every integer $k \geq r$.
  To this end, we compare the support functions of these polytopes on integer vectors.

  Let $\sigma \in \Z^{m+1}$, whose coordinates are indexed by $[0, m]$.
  For $I \subseteq [a]$ and $J \subseteq [b]$ with $\abs{I} = \abs{J} \geq 1$,
  let $B_{I, J} \coloneqq \sum_{i=0}^m t^{\sigma_i} A_i[I, J] x_i$, which is a linear symbolic matrix over $\F(t)$.
  For every positive integer $k$,
  the polynomial $\det B_{I, J}^{(k)}$ is obtained from $\det A[I, J]^{(k)}$ by substituting $t^{\sigma_i} X_i^{(k)}$ for $X_i^{(k)}$ for each $i \in [0, m]$,
  by which the monomial with an exponent vector $\alpha$ is multiplied by $t^{k \sigma^\top \proj_{k}(\alpha)}$.
  Hence, as observed in~\cite[(3.14)]{Hirai2025},
  \begin{align*}
    \frac{1}{k} \deg \det B_{I, J}^{(k)} = \max \Set*{\sigma^\top \proj_{k}(\alpha)}{\alpha \in \supp \det A[I, J]^{(k)}},
  \end{align*}
  where the maximum of the empty set is $-\infty$.
  For $I = J = \emptyset$, we let $\deg \det B_{\emptyset, \emptyset}^{(k)} \coloneqq 0$ and $\deg \Det B_{\emptyset, \emptyset} \coloneqq 0$;
  then the above equality also holds.
  Therefore, by the definitions of $\calP^{k}(A)$ and $\calQ(A)$ and by \Cref{prop:deg-det}~(3), we have
  \begin{align}\label{eq:support-functions}
    \max_{v \in \calP^{k}(A)} \sigma^\top v = \max_{I, J} \frac{1}{k} \deg \det B_{I, J}^{(k)},
    \qquad
    \max_{v \in \calQ(A)} \sigma^\top v = \max_{I, J} \deg \Det B_{I, J}
  \end{align}
  for every positive integer $k$,
  where the maxima are taken over the same pairs $(I, J)$ as in the definition of $\calP^{k}(A)$;
  see also~\cite[Theorem~3.17~(1)]{Hirai2025} for the latter.

  Now let $k \geq r$ be an integer.
  For every pair $(I, J)$ with $\abs{I} = \abs{J} \geq 1$ and every $i \in [m]$,
  the rank of the coefficient $t^{\sigma_i} A_i[I, J]$ of $B_{I, J}$ over $\F(t)$ is at most $\rank_\F A_i \leq r$.
  Hence $\deg \det B_{I, J}^{(k)} = k \deg \Det B_{I, J}$ by \Cref{prop:deg-det-attainment};
  this also holds for $I = J = \emptyset$.
  Therefore, by~\eqref{eq:support-functions}, we have $\max_{v \in \calQ(A)} \sigma^\top v = \max_{v \in \calP^{k}(A)} \sigma^\top v$ for every $\sigma \in \Z^{m+1}$.
  Since $\calQ(A)$ and $\calP^{k}(A)$ are polytopes, this implies $\calQ(A) \subseteq \calP^{k}(A)$ by~\cite[Lemma~2.2]{Hirai2026}.
\end{proof}

\section*{Acknowledgments}
This work was supported by JSPS KAKENHI Grant Numbers
\mbox{JP22K17854}, \mbox{JP24K21315}, and \mbox{JP24K02901}, Japan.

\section*{Declaration on the use of generative AI tools}
The author used ChatGPT-6.1 Sol and Claude Fable 5.2 for discussions, drafting, and proofreading.
The author has carefully reviewed and verified the entire manuscript and takes full responsibility for all of its contents.

\bibliography{ref}

\begin{thebibliography}{10}

\bibitem{Amitsur1966}
S.~A. Amitsur.
\newblock Rational identities and applications to algebra and geometry.
\newblock {\em Journal of Algebra}, 3(3):304--359, 1966.

\bibitem{Cohn1995}
P.~M. Cohn.
\newblock {\em Skew Fields: Theory of General Division Rings}, volume~57 of
  {\em Encyclopedia of Mathematics and its Applications}.
\newblock Cambridge University Press, Cambridge, 1995.

\bibitem{Derksen2017}
H.~Derksen and V.~Makam.
\newblock Polynomial degree bounds for matrix semi-invariants.
\newblock {\em Advances in Mathematics}, 310:44--63, 2017.

\bibitem{Derksen2018}
H.~Derksen and V.~Makam.
\newblock On non-commutative rank and tensor rank.
\newblock {\em Linear and Multilinear Algebra}, 66(6):1069--1084, 2018.

\bibitem{Dieudonne1943}
J.~Dieudonn{\'e}.
\newblock Les d{\'e}terminants sur un corps non commutatif.
\newblock {\em Bulletin de la Soci{\'e}t{\'e} Math{\'e}matique de France},
  71:27--45, 1943.

\bibitem{Edmonds1967}
J.~Edmonds.
\newblock Systems of distinct representatives and linear algebra.
\newblock {\em Journal of Research of the National Bureau of Standards},
  71B(4):241--245, 1967.

\bibitem{Fortin2004}
M.~Fortin and C.~Reutenauer.
\newblock Commutative/noncommutative rank of linear matrices and subspaces of
  matrices of low rank.
\newblock {\em S{\'e}minaire Lotharingien de Combinatoire}, 52:B52f, 2004.

\bibitem{Garg2020}
A.~Garg, L.~Gurvits, R.~Oliveira, and A.~Wigderson.
\newblock Operator scaling: theory and applications.
\newblock {\em Foundations of Computational Mathematics}, 20:223--290, 2020.

\bibitem{Gille2006}
P.~Gille and T.~Szamuely.
\newblock {\em Central Simple Algebras and {G}alois Cohomology}.
\newblock Cambridge University Press, Cambridge, 2006.

\bibitem{Gurvits2004}
L.~Gurvits.
\newblock Classical complexity and quantum entanglement.
\newblock {\em Journal of Computer and System Sciences}, 69(3):448--484, 2004.

\bibitem{Hamada2021}
M.~Hamada and H.~Hirai.
\newblock Computing the nc-rank via discrete convex optimization on {CAT}(0)
  spaces.
\newblock {\em SIAM Journal on Applied Algebra and Geometry}, 5(3):455--478,
  2021.

\bibitem{Hirai2019}
H.~Hirai.
\newblock Computing the degree of determinants via discrete convex optimization
  on {E}uclidean buildings.
\newblock {\em SIAM Journal on Applied Algebra and Geometry}, 3(3):523--557,
  2019.

\bibitem{Hirai2026}
H.~Hirai.
\newblock On integral polytopes related to {E}dmonds' problem.
\newblock arXiv:2609.19703v1, 2026.

\bibitem{Hirai2025}
H.~Hirai, Y.~Iwamasa, T.~Oki, and T.~Soma.
\newblock Algebraic combinatorial optimization on the degree of determinants of
  noncommutative symbolic matrices.
\newblock {\em Mathematical Programming}, 213:941--984, 2025.

\bibitem{Hrubes2015}
P.~Hrube{\v{s}} and A.~Wigderson.
\newblock Non-commutative arithmetic circuits with division.
\newblock {\em Theory of Computing}, 11(14):357--393, 2015.

\bibitem{Huszar2021}
A.~Huszar.
\newblock Non-commutative rank and semi-stability of quiver representations.
\newblock arXiv:2111.00039v1, 2021.

\bibitem{Ivanyos2015}
G.~Ivanyos, M.~Karpinski, Y.~Qiao, and M.~Santha.
\newblock Generalized {W}ong sequences and their applications to {E}dmonds'
  problems.
\newblock {\em Journal of Computer and System Sciences}, 81(7):1373--1386,
  2015.

\bibitem{Ivanyos2010}
G.~Ivanyos, M.~Karpinski, and N.~Saxena.
\newblock Deterministic polynomial time algorithms for matrix completion
  problems.
\newblock {\em SIAM Journal on Computing}, 39(8):3736--3751, 2010.

\bibitem{Ivanyos2017}
G.~Ivanyos, Y.~Qiao, and K.~V. Subrahmanyam.
\newblock Non-commutative {E}dmonds' problem and matrix semi-invariants.
\newblock {\em Computational Complexity}, 26:717--763, 2017.

\bibitem{Ivanyos2018}
G.~Ivanyos, Y.~Qiao, and K.~V. Subrahmanyam.
\newblock Constructive non-commutative rank computation is in deterministic
  polynomial time.
\newblock {\em Computational Complexity}, 27:561--593, 2018.

\bibitem{Kabanets2004}
V.~Kabanets and R.~Impagliazzo.
\newblock Derandomizing polynomial identity tests means proving circuit lower
  bounds.
\newblock {\em Computational Complexity}, 13:1--46, 2004.

\bibitem{Lovasz1989}
L.~Lov{\'a}sz.
\newblock Singular spaces of matrices and their application in combinatorics.
\newblock {\em Boletim da Sociedade Brasileira de Matem{\'a}tica}, 20:87--99,
  1989.

\bibitem{Oki2025}
T.~Oki and T.~Soma.
\newblock Algebraic algorithms for fractional linear matroid parity via
  noncommutative rank.
\newblock {\em SIAM Journal on Computing}, 54(1):134--162, 2025.

\bibitem{Schwartz1980}
J.~T. Schwartz.
\newblock Fast probabilistic algorithms for verification of polynomial
  identities.
\newblock {\em Journal of the ACM}, 27(4):701--717, 1980.

\bibitem{Shoup1994}
V.~Shoup.
\newblock Fast construction of irreducible polynomials over finite fields.
\newblock {\em Journal of Symbolic Computation}, 17(5):371--391, 1994.

\bibitem{Stacks}
{The Stacks Project Authors}.
\newblock The {S}tacks project.
\newblock \url{https://stacks.math.columbia.edu}, 2026.

\bibitem{Gathen2013}
J.~von~zur Gathen and J.~Gerhard.
\newblock {\em Modern Computer Algebra}.
\newblock Cambridge University Press, Cambridge, third edition, 2013.

\bibitem{Zippel1979}
R.~Zippel.
\newblock Probabilistic algorithms for sparse polynomials.
\newblock In E.~W. Ng, editor, {\em Symbolic and Algebraic Computation},
  volume~72 of {\em Lecture Notes in Computer Science}, pages 216--226, Berlin,
  Heidelberg, 1979. Springer.

\end{thebibliography}

\end{document}